\documentclass{amsart} 
\usepackage{latexsym,amssymb,graphicx}
       
\newtheorem{thm}{Theorem}[section]
\newtheorem{lem}[thm]{Lemma}
\newtheorem{prop}[thm]{Proposition}
\newtheorem{cor}[thm]{Corollary} 
\newtheorem{de}[thm]{Definition} 
\newtheorem{rem}[thm]{Remark}
\newtheorem{ex}[thm]{Example}

\newcommand{\BZ}{{\mathbb{Z}}}
\newcommand{\BQ}{{\mathbb{Q}}}
\newcommand{\BO}{{\mathcal{O}}}
\newcommand{\BD}{{\mathcal{D}}}

\newcommand{\BB}{{\mathcal{B}}}

\newcommand{\BG}{{\mathcal{G}}}

\newcommand{\BJ}{{\mathcal{J}}}
\newcommand{\bG}{{\mathbb{G}}}

\newcommand{\BF}{{\mathbb{F}}}

\newcommand{\BS}{{\mathcal{S}_p}}
\newcommand{\BSplus}{{\mathcal{S}^+_p}}

\newcommand{\BSs}{{{\mathcal{S}}_p^{\sharp}}}
\newcommand{\BSps}{{{\mathcal{S}}_p^{+\sharp}}}

\newcommand{\Vp}{{V_p}}

\newcommand{\Si}{{\Sigma}}
\newcommand{\oi}{{\mathfrak{o}}}
\DeclareMathOperator{\rad}{rad}

\newcommand{\bb}{{\mathfrak{b}}}
\newcommand{\bg}{{\mathfrak{g}}}

\newcommand{\bB}{{\mathbb{B}}}
\newcommand{\Bs}{{{\mathcal{S}}}}
\newcommand{\Bss}{{{\mathcal{S}}^{\sharp}}}
\newcommand{\bo}{{\mathfrak{o}}} 
\newcommand{\bc}{{s}}

\def\mapright#1{\smash{\mathop{\longrightarrow}\limits^{#1}}}

\begin{document}
\title{Integral Lattices in TQFT}

\author{ Patrick M. Gilmer}
\address{Department of Mathematics\\
Louisiana State University\\
Baton Rouge, LA 70803\\
USA}
\email{gilmer@math.lsu.edu}
\thanks{The first author was partially supported by NSF-DMS-0203486}
\urladdr{www.math.lsu.edu/\textasciitilde gilmer/}

\author{ Gregor Masbaum}
\address{Institut de Math{\'e}matiques de Jussieu (UMR 7586 du CNRS)\\
Universit{\'e} Paris 7 (Denis Diderot) \\
Case 7012\\
2, place Jussieu\\
75251 Paris Cedex 05\\
FRANCE }
\email{masbaum@math.jussieu.fr}
\urladdr{www.math.jussieu.fr/\textasciitilde masbaum/}

\date{
February 26, 2007 
\hspace{0.5cm} First edition: October 31, 2004.}

\begin{abstract} We find explicit bases for 
naturally defined  lattices over a ring of algebraic integers in
the $SO(3)$-TQFT-modules of surfaces 
at roots of unity of odd prime order.  
Some applications relating quantum invariants to classical
$3$-manifold topology are given.
\end{abstract}

\maketitle
\tableofcontents

\section{\   Introduction} \label{sec.intro}

 Let $p$ be an odd prime. Based on integrality results for the quantum
 $SO(3)$-invariant of closed oriented $3$-manifolds 
      \cite{Mu2,MR}, an
 integral TQFT-functor $\BS$ was defined in \cite{G} and
 \cite{GMW}. It  associates to a closed surface $\Si$ a free
lattice\footnote{A lattice over a Dedekind domain is a finitely generated
torsion-free module. 
In general, a lattice need not be  free,
and the freeness of the lattices $\BS(\Si)$ is a non-trivial fact.} 
 $\BS(\Si)$ over the cyclotomic ring 
\begin{equation*}\label{cyclo}
\BO=\left\{
\begin{array}{cl}
\BZ[\zeta_p] &\text{if $p \equiv -1 \pmod{4}~,$} \\
\BZ[\zeta_{p},i]=\BZ[\zeta_{4p}] &\text{if $p \equiv 1 \pmod{4}~,$}\\
\end{array}\right.
\end{equation*}
where $\zeta_n$ is a primitive $n$-th root of unity.
 The lattice $\BS(\Si)$  carries an $\BO$-valued non-degenerate
 hermitian form which will be denoted by $(\ ,\ )_\Si$. Here,
 non-degenerate means that the form induces an injective adjoint map
 $\BS(\Si)\rightarrow \BS(\Si)^*$. The lattice $\BS(\Si)$ also carries
 a linear  representation of an 
appropriately extended
mapping class group of $\Si$; moreover, this representation preserves the hermitian form on $\BS(\Si)$.

The integral lattices $\BS(\Si)$ have a natural definition in terms of the
vector-valued quantum $SO(3)$-invariants for $3$-manifolds with boundary (see
below). Their existence thus reflects interesting
structural properties of quantum invariants and TQFT's. 
 The main aim of the present paper is to give an explicit description of
 these lattices for surfaces of arbitrary genus, 
possibly disconnected,
by describing bases for them.

Here, we allow surfaces to be equipped with a (possibly empty) collection of 
colored banded points (a banded point is an 
embedded 
oriented arc), where a color is an integer $i\in\{0,1,\ldots,  p-2\}$.
But we require the sum of colors on the colored points of any component 
of a surface to be even (this is a feature of the $SO(3)$-theory we
consider). Before, explicit bases of $\BS(\Si)$ were known
only in the case where the surface $\Si$ 
is connected and 
 has genus one or two with no colored points \cite{GMW}.

The bases we find display some nice ``graph-like'' structure which we
believe should generalize to other TQFT's, at least
for those associated to integral modular categories as defined in \cite{MW}.
 One can ask 
whether the graph-like bases we have found might be related to canonical  bases in representation theory.

An interesting new feature is that the usual tensor
product axiom of TQFT holds only with some modification for the 
lattices $\BS(\Si)$: it turns out that the lattice associated to a
disjoint union of surfaces is sometimes bigger than the tensor product
of the lattices associated to the individual components, although this
phenomenon does not happen for surfaces without colored points.

 We remark that Kerler \cite{K} has announced an integral
 version
of the $SO(3)$ TQFT at $p=5$, but details  have not yet appeared. 
 Chen and Le have recently constructed integral TQFTs  from quantum 
groups \cite{CL, C}, although without describing explicit bases for the
modules associated to surfaces.

We give two topological applications of our theory.
In Section \ref{sec.cut}, we prove a conjecture which relates the cut number of a $3$-manifold $M$ to 
 divisibility properties of its quantum invariants. The cut number is  the same as the co-rank of $\pi_1(M)$. Thus this result relates quantum invariants to the fundamental group of a $3$-manifold.
   
Our second application concerns the
Frohman and  Kania-Bartoszynska ideal of  $3$-manifolds with boundary  \cite{FK} which can be used to show that one $3$-manifold does not embed in another.  This ideal is hard to compute directly from its definition, except in very special circumstances, as its definition involves the quantum invariants of infinitely many $3$-manifolds.
However our integral bases allow us to give an explicit finite set of generators for this ideal (see Section \ref{sec.FKB}). We apply this method to exhibit  a family of examples of 
$3$-manifolds with the homology of a solid torus which cannot embed in $S^3.$

The integral TQFT-functor $\BS$ is a refinement of the
$SO(3)$-TQFT-functor $V_p$ constructed in \cite{BHMV2}, which in turn
is, in some sense, an alternative version of a special case of the
Reshetikhin-Turaev theory \cite{RT} of quantum invariants of $3$-manifolds. In
particular, the lattice $\BS(\Si)$ is defined in \cite{G,GMW} as an
$\BO$-submodule of the TQFT-vector space $V_p(\Si)$. The latter is
defined over the quotient field of $\BO$ and exists also when $p$ is
not prime. In fact, $V_p(\Si)$ can be defined (and is a free module) already over the ring $\BO[\frac 1 p]$, see \cite{BHMV2},
and this is the version of $V_p$ that we refer to in the rest of the paper.

As a submodule of  $V_p(\Si)$, the lattice $\BS(\Si)$
is simply defined as the $\BO$-span of the
vectors associated to $3$-manifolds $M$ with boundary $\partial M
=\Si$, with the condition that $M$ should have no closed
components. It is clear from this definition that the extended mapping
class group acts on $\BS(\Si)$. The fact that $\BS(\Si)$ is a free
lattice of rank equal to the dimension of $V_p(\Si)$ is shown in
\cite{G}. The hermitian form $(\ ,\ )_\Si$ is obtained by rescaling
the natural form on $V_p(\Si)$ given by gluing manifolds together
along their boundary and computing the invariant of the closed
manifold so obtained. That such a rescaling is possible depends
crucially on the integrality result for this invariant of closed
$3$-manifolds due 
to  H.~Murakami \cite{Mu2} and Masbaum-Roberts
\cite{MR} (see \cite{G} for more details). 
\footnote{Using deep results of Habiro, integrality of the $SO(3)$-invariant for all $3$-manifolds
  and any odd $p$ was proven in 2006  by Beliakova-Le
  \cite{BL} (see also Le \cite{L}).}

Bases of the free $\BO[\frac 1 p]$-module $V_p(\Si)$ are well
understood in terms of admissible colorings of uni-trivalent graphs. Here,
any uni-trivalent graph which is the spine of a handlebody with boundary 
 $\Si$, and with the univalent vertices meeting $\Si$ in
the colored points, may be used. The $\BO$-span of such a
{\em graph basis} is a sublattice of $\BS(\Si)$, but this sublattice 
is almost never invariant under the mapping class
group, and hence cannot be equal to the whole integral lattice 
$\BS(\Si)$. One might hope that a basis of $\BS(\Si)$ could be
obtained by rescaling the elements of a graph basis in some way, but 
this is not
the case. Still, the situation is actually rather nice. We will show that
the lattice $\BS(\Si)$ admits what we call {\em graph-like} bases
associated to a special kind of uni-trivalent graph which we call a 
{\em lollipop tree}. Roughly speaking, a graph-like basis is obtained
from the usual graph basis associated to the lollipop tree by taking certain linear combinations, 
followed by some overall rescaling depending on the colors. The nice
thing is that the linear
combinations are taken independently in each handle. For precise
definitions and a statement of the result in the case of connected 
surfaces, see
Section~\ref{sec.thm}.

We remark here that 
for connected surfaces, 
$\BS(\Si)$ has a simple skein-theoretical description shown in \cite{GMW}, namely as the $\BO$-span 
in  $V_p(\Si)$
 of banded links and graphs in a handlebody colored in a certain way.
 This description will be given below in 
Proposition \ref{vgraph}.
For the purpose of the present paper, this description  
can
be taken as the definition of $\BS(\Si)$. 
Moreover, our method of constructing a basis will give an independent proof that $\BS(\Si)$ is indeed a free lattice.

For disconnected surfaces, we describe a basis of $\BS(\Si)$ in Section
\ref{discon.sec}, 
where we also discuss the modified tensor product axiom.

In the case of surfaces of genus one and two without colored
points, the natural hermitian form $(\ ,\ )_\Si$ was shown in
\cite{GMW} to be
unimodular (here, unimodular means that the adjoint map is not only 
injective but  is an isomorphism). 
 As already observed there, this property no longer holds in higher
 genus. 
It is then natural to consider the dual lattice $\BSs(\Si)$, defined as 
\[ \BSs(\Si) =\{x\in 
V_p(\Si)\ 
|\ (x,y)_\Si\in \BO \text{ for
 all } y\in \BS(\Si)\}~. \]  
Note that $\BS(\Si)\subset  \BSs(\Si)$, with equality if and only if the form is unimodular.\footnote{It is easy to check that since $\BS(\Si)$ is a 
free
lattice, so is $\BSs(\Si)$. In fact, $\BSs(\Si)$ is isomorphic 
as an $\BO$-module to the dual of 
the conjugate of 
$\BS(\Si)$, justifying the terminology.}

The main result about the dual lattices $ \BSs(\Si)$ is that they also admit
graph-like bases. In fact, such bases can be described as rescalings of
graph-like bases for $\BS(\Si)$ (see again Section~\ref{sec.thm} for a
precise statement). Moreover, the dual lattices play an important
role in the proof that our bases are indeed bases, 
as the proof 
proceeds by studying the two lattices simultaneously.

The dual lattice $\BSs(\Si)$ is, of course, also preserved by the
action of the mapping class group. Note that when one expresses the
mapping class group representation on $ \BS(\Si)$ or on  $\BSs(\Si)$ in our integral bases, all matrix
coefficients will be algebraic integers.

 One also obtains a representation 
 by isometries 
 (of the associated torsion ``linking form'')
   on the quotient $\BSs(\Si)/\BS(\Si)$, which admits a simple
   description,
at least if $p\equiv -1\pmod{4}$: 
  When $\Si$ is connected and has no colored points,
  $\BSs(\Si)/\BS(\Si)$ 
is a skew-symmetric inner product space over the finite field $\BF_p$.
If $p\equiv 1\pmod{4}$, a similar statement holds with a little more
effort for a refined theory, see 
Section~\ref{torsion.sec}
for details.
It should be mentioned, however, that here we merely begin the study of these
representations of mapping class groups on torsion modules; we hope to
return to this matter elsewhere.

 {\em Acknowledgments.} We would like to thank the following for hospitality and support:
the visiting experts program in mathematics of the Louisiana Board of Regents (LEQSF(2002-04)-ENH-TR-13), 
Knots in Poland 2003 at The Banach Center and at the Mathematical Research and Conference Center in
Bedlewo, the Mini-Workshop: Quantum Topology in Dimension Three at
Mathematisches Forschungsinstitut Oberwolfach,  and also the
Mathematics Institute of Aarhus University and the MaPhySto network of
the Danish National Research Foundation.

{\em Note added March 2006.} As a further application of the integral lattices $\BS(\Si)$, it will be shown in \cite{M?} that when one expresses the mapping class group representations  in our graph-like bases associated to lollipop trees, then the representations (at least when restricted to the Torelli group) have a perturbative limit as $p\rightarrow \infty$, in much the same sense as Ohtsuki's power series invariant \cite{Oh} of homology spheres is a limit of the quantum $SO(3)$-invariants at roots of prime order.

 {\em Notations and  Conventions.}
Throughout the paper,
$p\geq 5$  
will be a
prime integer, and we
put $d=(p-1)/2$. 
\footnote{We assume $p\geq 5$ because the color $2$, which will play an important
role in our construction, does not exist for $p=3$. But the theory 
for $p=3$ is trivial anyway \cite{BHMV2}.}
Recall that our ground ring $\BO$ contains a
primitive $p$-th root of unity $\zeta_p$. We define $h\in \BO$ by  
\begin{equation} h=1-\zeta_p~. \end{equation} One has that $h^{p-1}$ is
a unit times $p$, so that $\BO[\frac 1 p]=\BO[\frac 1 h]$. 
For skein
theory purposes, we put $A= -
\zeta_p^{d+1}$. This is a primitive $2p$-th root of
unity such that $A^2=\zeta_p$. The quantum integer $[n]$ is defined  by $[n]=(\zeta_p^{n}-\zeta_p^{-n})/(\zeta_p-\zeta_p^{-1})$. 
We also fix $\BD\in\BO$ such that 
\begin{equation} \label{defD} \BD^2 = \frac{-p }{(\zeta_p-\zeta_p^{-1})^2}~. 
\end{equation} One has that 
$\mathcal{D}$ is a unit times $h^{d-1}.$

Unless otherwise stated all manifolds that we consider are assumed to be compact and oriented.
If a surface $\Si$ is equipped with a particular collection of colored
 banded points, we denote the latter by $\ell(\Si).$ 
When writing $V_p(\Si)$, $\BS(\Si)$, {\em etc.}, we let $\Si$ stand for the surface equipped with its given colored
points.

As usual in TQFT, our surfaces and $3$-manifolds are equipped with 
an additional structure to resolve the ``framing anomaly''. As it is 
well-known 
how to do this \cite{BHMV2,W,T}, and the additional structure is 
basically irrelevant for  integrality questions, we postpone further 
discussion of
this ``framing'' issue until Section~\ref{mcg.sec} where some details
of the construction will be needed.

\section{\   Skein theory and the definition of $\BS(\Si)$}\label{skein.sec1}

Unless otherwise stated, 
by skein module we will mean the Kauffman Bracket skein module over $\BO[\frac 1 p]=\BO[\frac 1 h]$.
Recall that this skein module  of a $3$-manifold $M$ 
is  the free $\BO[\frac 1 h]$  module on the banded links in $M$
modulo the well-known Kauffman relations and isotopy.  We denote this
module by 
$K(M).$  Elements   of  $K(M)$ can be described by colored trivalent
banded graphs, which can be expanded to linear combinations of banded
links in the familiar way.  In particular, a strand colored, say $a$,
is replaced by $a$ parallel strands with the Jones-Wenzl idempotent,
denoted $f_a$, inserted. This is a particular linear combination of $a$-$a$ tangles, and is sometimes denoted as a rectangle with $a$ inputs on each of the long sides. See for instance  \cite{BHMV2}.
We also need to consider the relative
skein module of a $3$-manifold whose boundary is equipped with some
banded colored points $\ell(\partial M)$. Here one takes the free
$\BO[\frac 1 h]$ module on those links 
(or rather, tangles)
which are expansions of colored
graphs which meet the boundary nicely in 
the colored points.
In this case
we only use isotopy relative the boundary in the relations. This
relative skein module is denoted $K(M, \ell(\partial M))$.
An element of this module is often denoted $(M,L)$ where $L$ stands
for the colored link or graph in $M$.

Suppose the surface 
$\Si$ is
equipped with a (possibly empty) collection of colored points
$\ell(\Si)$.  
The $\BO[\frac 1 p]$-module $V_p(\Si)$ can be described as follows. 
For any 
$3$-manifold $M$ with boundary $\Si$,
there is a surjective
map from $K(M,\ell(\Si))$ to $V_p(\Si)$.~\footnote{There is no 
connectivity hypothesis on $M$ since $V_p$ satisfies the tensor 
product axiom. Here we use our assumption that the sum of the
  colors of the banded points on every connected component of every
  surface is even.}
The image of the skein class represented by $(M,L)$ is denoted
$[(M,L)]$. Here we think of $(M,L)$ as a cobordism from $\emptyset$ to
$\Si.$ 
Suppose $M'$ is a second 
$3$-manifold with boundary  $\Si$, then
a sesquilinear form 
 $\langle \ ,\ \rangle_{M, M'}: K(M,\ell(\Si)) \times K(M',\ell(\Si)) \rightarrow \BO[\frac 1 p]$  is defined by \[  \langle(M,L), (M',L')\rangle_{M, M'} =  
 \langle (M \cup_\Si -M', L \cup_{\ell(\Si)} -L')\rangle_p
  \]   Here 
$\langle \ \rangle_p$ denotes the quantum invariant of a
 closed $3$-manifold,
and the minus sign indicates reversal of orientation.
 The kernel of the map  $K(M,\ell(\Si)) \rightarrow V_p(\Si)$ is the left 
radical of the form $\langle \ ,\ \rangle_{M, M'}$. Moreover, 
this form 
induces the canonical 
nonsingular  hermitian  form
 $\langle \ ,\ \rangle_{\Si}:V_p(\Si) \times V_p(\Si) \rightarrow
 \BO[\frac 1 p]$ 
(which is independent of $M$ and $M'$).
 All the results of this paragraph appear in
 \cite{BHMV2}.

\begin{de}\label{defS}{\em Given a closed surface $\Si$ with possibly a collection
 of colored banded points $\ell(\Si)$, we define $\BS(\Si)$ to be the
 $\BO$-submodule of $V_p(\Si)$ generated by all vectors $[(M,L)]$
 where  
$M$ is any $3$-manifold with boundary $\Si$ having no closed connected components, and the colored graph $L\subset M$ meets $\Si$ nicely in $\ell(\Si)$.
}\end{de}

As shown in \cite{GMW}, $\BS(\Si)$ has a skein-theoretical
description, as follows.  First, we let $v$ denote the skein element in $K(S^1 \times D^2)$ described by
$h^{-1} ( 2+ z)$ where $z$ is $S^1 \times \{0\}$ with 
standard framing.
Thus $v$ denotes a skein class and also the element
in 
$V_p(
S^1 \times S^1)$ which 
it
represents, depending on context. 
(The $v$ used in \cite{GMW} differs by a unit from  the $v$ used
here.)
Coloring a link component $v$ is shorthand for replacing this component by the linear combination $h^{-1} ( 2+ z)$  and expanding linearly.

Next, by a {\em $v$-graph} in  $3$-manifold with possibly some colored
points in the boundary, we will mean a banded colored graph in $M$
which agrees with $\ell(\partial M)$ on the boundary together with
possibly some other 
banded
link components which have been colored with
$v$. By \cite[Prop. 5.6 and Cor. 7.5 ]{GMW},  we have that:

\begin{prop} \label{vgraph} Suppose $\Si$ is a connected surface. 
Choose a connected  $3$-manifold $M$ 
with boundary $\Si$.  Then $\BS(\Si)$ is generated over $\BO$ by elements  represented by $v$-graphs in $M$ which meet the boundary in the colored points of $\Si$. 
\end{prop}

For connected surfaces, the above can be taken as an alternative 
definition of $\BS(\Si)$. 
There is also a version of this for surfaces which are not connected
but we delay stating it until  section \ref{discon.sec} when it is
needed.

\section{\   Lollipop trees and the small graph basis of $V_p(\Si)$}

We let $\Si$ denote the boundary of a
genus $g$ handlebody $H_g$ and fix a particular collection of colored
banded points in $\Si$ which we denote by $\ell(\Si).$
 A basis of $V_p(\Si)$ can be described by 
the $p$-admissible 
(see below)
colorings of any 
uni-trivalent 
banded graph $G$ having the same homotopy type as 
the handlebody $H_g$
and which meets the boundary at $\ell(\Si)$
in the univalent vertices of $G$. 
Moreover,  the colorings
should extend the given colorings at these banded colored points 
 \cite{BHMV2}.  
 As usual, a colored graph with an edge colored zero is identified with the same graph without that edge (similarly for zero-colored points in $\ell(\Si)$).

\begin{figure}[h]
\includegraphics[width=2in]{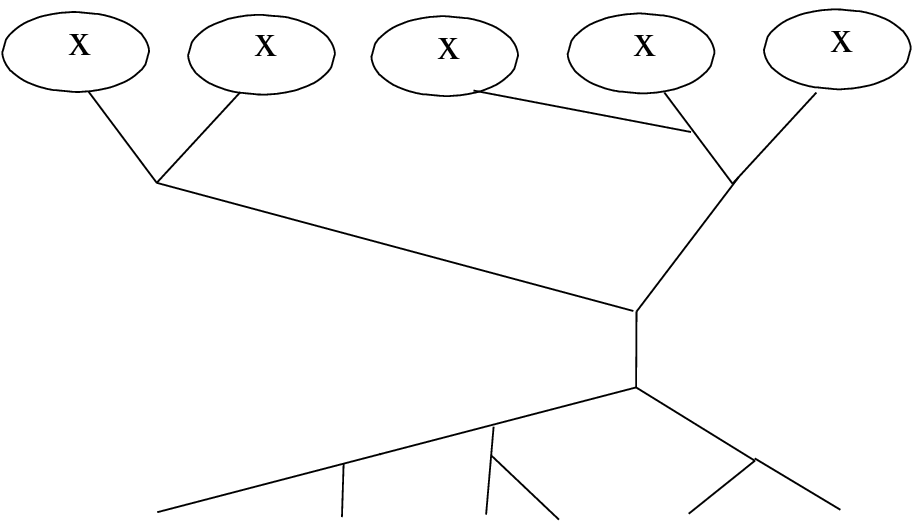}
\caption{} \label{t}
\end{figure}

 For example, we can take the graph  in  Figure \ref{t} when $g=5$ and $\ell(\Si)$ has 6 points.
The banding of this graph lies  in the plane. 
The points $\ell(\Si)$ are depicted
at the bottom of the diagram. The x's in this and later diagrams denote  holes in $H_g$.

A $p$-admissible coloring of $G$ is an assignment of colors to the edges of $G$ such that at every vertex of $G$, the three colors $i,j,k$, say, meeting at that vertex satisfy the conditions 
\begin{eqnarray*}\label{adm1} i+j+k &\equiv& 0 \pmod 2\\
\label{adm2} |i-j|\ \leq& k &  \leq \ i+j\\
\label{adm3} i+j+k &\leq& 2p-4=4d-2
\end{eqnarray*}

To a $p$-admissible coloring of $G$ one associates in the usual way a skein element in the handlebody $H_g$, by replacing the edges of $G$ with appropriate Jones-Wenzl idempotents. Identifying now the boundary of the handlebody with our surface $\Si$, this skein element defines in turn a vector in $V_p(\Si)$.
The vectors associated to $p$-admissible colorings where, in addition,
the colors 
satisfy a parity condition, form a basis of $V_p(\Si)$
\cite[Theorem 4.14]{BHMV2}. \footnote{In the exceptional case where $\Si$ is a two-sphere, with
 only one 
even-colored
point,  there is no such graph. Thus
 $V_p(\Si)$ is zero if the color is nonzero, and $\BO[\frac 1 p]$ if the color is zero.}

In Proposition~\ref{2.1} below, we will describe a different basis of $V_p(\Si)$, where the parity
condition is replaced by a ``smallness'' condition.  For this, we must
restrict our graph to be a lollipop tree, defined as follows. 

\begin{de}{\em Let $G$ be a uni-trivalent graph as above. Let $g$ be
    the first Betti number of $G$ and let 
$\bc$ be the number of its
    univalent vertices. (Of course, $g$ is the genus of $\Sigma$ and
    $\bc$ is the number of colored points in $\ell(\Si)$.) Then $G$ is
    called a {\em lollipop tree} provided it satisfies the following
    conditions.
\begin{itemize}
\item[(i)] $G$ has exactly $g$ loop edges, so that the complement of
  these loop edges in $G$ is a tree $T$.
\item[(ii)] If $\bc>0$, there must be a single edge of the tree $T$ 
called 
the {\em trunk edge} (or simply the {\em trunk}) 
with the property that  if we remove the interior of the trunk from $T$,
we obtain the  disjoint union of two trees: one which meets every loop
edge  and one which 
contains every uni-valent vertex. (If $\bc=1$, this second tree consists of a single point.) 
\end{itemize}
}\end{de}

Note that a lollipop tree is not actually a tree. 
We chose to call it so because of  the special case where there is only one loop edge and the tree $T$ consists of just one edge; in this case the graph $G$ looks somewhat like a lollipop. 
  
For example, the 
graph
of Figure \ref{t} is a lollipop tree.

\begin{prop}\label{2.1} The vectors associated to $p$-admissible
  colorings of a lollipop tree $G$, where the loop edges are assigned colors in the 
interval $[0,d-1]$, form a basis of $V_p(\Si)$.
\end{prop}

We will refer to this basis as the {\em small graph basis} and denote it by $\BG$.  Here, the adjective ``small'' refers to the colors on the loop edges.

\begin{proof} 
By \cite[Theorem 4.14]{BHMV2}, we have a basis by taking all
$p$-admissible colorings with even colors on the loop edges. (Observe
that the parity of the colors of the edges of the sub-graph
$T$ is imposed by the colors of $\ell(\Si)$.) Lemma 8.2 of \cite{GMW}
shows how to replace even colors on the loop edges by small ones,
{\em i.e.,} colors in $[0,d-1]$. Specifically, the lemma dealt with
the case $g=2,\bc=0$, but the same argument works in general. \end{proof}

\section{\   Bases for $\BS(\Si)$ and $\BSs(\Si)$ for connected surfaces}\label{sec.thm}

We are now ready to state our results concerning graph-like bases of
$\BS(\Si)$ and $\BSs(\Si)$. Fix a lollipop tree $G$ and define $g$ and
$\bc$  as in the preceding
section.

We will use the following labelling of the edges of $G$. Recall that
the complement of the loop edges of $G$ is a tree $T$. We call an edge
of $T$   {\em a stick edge}  if it is incident with a loop edge of $G$. 
An edge is called {\em ordinary} otherwise. 
We denote their colors by $2a_1, \ldots, 2a_g$ for the stick edges
(which will always have an even color),  and by 
$c_1,c_2, \ldots $
for
the ordinary edges. Here
 \begin{equation} \label{smallbasis1} 0\leq a_i \leq d-1~.
\end{equation}
The case $g=2,\bc=0$ is special as there is only one stick edge. In this
case we put $a_1=a_2$ but both $a_1$ and $a_2$ are to be entered in  
Eqs.~(\ref{basedef0}), 
(\ref{basedef}), and (\ref{basedefsharp}) below.
  
The color of the trunk is always
even; we denote it by $2e$.  If $\bc=0$,
we also set $e$ to be zero.   We note that the trunk is usually an
ordinary edge but it is a stick edge when $g=1$ and 
$\bc\geq 1$.

We denote the colors of the loop edges by $a_i+b_i$ ($i=1,\ldots,
g$). Here, the loop edge colored $a_i+b_i$ is incident to the 
stick
edge of $T$ labelled $2a_i$. 
Note $b_i\geq 0$ by $p$-admissibility. Moreover, since loop edges should have small colors, we have
 \begin{equation} \label{smallbasis2} 0\leq b_i\leq d-1-a_i~.
\end{equation}
The elements of the small graph basis $\BG$ will be denoted by
$\bg(a,b,c)$, where $a=(a_1,\ldots, a_g)$, $b=(b_1,\ldots, b_g)$, and
 $c=(c_1,c_2, \ldots )$. The index set is precisely the set of 
$(a,b,c)$ satisfying conditions \eqref{smallbasis1} and \eqref{smallbasis2}, and such that  $(2a,c)$ is a $p$-admissible coloring of the tree $T$
extending the given coloring of $\ell(\Si)$.
 We will refer to $(a,b,c)$ as  a small coloring
of $G$.
\vskip 8pt

{\em The basis of $\BS(\Si)$.} Our basis of $\BS(\Si)$ will be denoted by $\BB$. It consists of vectors $\bb(a,b,c)$ indexed by the same set as the small graph basis $\BG$.
  They are defined as follows.  Recall that  $h=1-\zeta_p$.  
For $x\in \mathbb R$, we use the notation $\lfloor x
\rfloor$ (resp. $\lceil x \rceil$) to denote the greatest integer
$\leq x$ (resp. the smallest integer $\geq x$).

If $b=0$, the vector $\bb(a,0,c)$ is just a rescaling of $\bg(a,0,c)$:  
\begin{equation} \label{basedef0} \bb(a,0,c) = 
h^{- \lfloor \frac 1 2 (-e+\sum_i a_i)
\rfloor } \bg(a,0,c) \end{equation} 
(Observe that one always has $e\leq \sum_i a_i$.) 
If $b>0$, the vector $\bb(a,b,c)$ 
is conveniently described using multiplicative notation, as follows.

Think of the handlebody $H_g$ as $P\times I$ where $P$ is a $g$-holed disk thought of as a regular neighborhood of a planar embedding of the banded graph $G$. This endows the 
absolute
 skein module of $H_g$ with an algebra structure, where multiplication is given by putting one skein element on top of the other. 
 Similarly, relative skein modules  of $H_g$ are modules over the absolute skein module.
 Note that this multiplication is non-commutative in general.  This also induces an algebra structure on $V_p(\Si)$ in the case $\ell(\Si)$ is empty.  Although these algebra/module structures are not canonically associated to the surface $\Si$, they are well-defined once 
 $\Si$ has been identified with the boundary of $P\times I$ 
 with all colored points, say, at level  one-half.  
 Similarly, we obtain  $\BO$-algebra/module structures on the lattices $\BS(\Si)$, which will be used throughout the paper. 

For $i=1,\ldots,g$, let $z_i$ denote the skein element represented by a  circle around the $i$-th hole of the $g$-holed disk $P$ and with framing parallel to $P$. We have  
\[ z_i=\bg((0,\ldots,0),(0,\ldots, 0,1,0,\ldots 0), (0,\ldots ))~. \]
(The only non-zero  coefficient sits at the $i$-th entry of the $g$-tuple $b$.) We define 
 \begin{equation*} v_i=h^{-1}(2+z_i) \end{equation*} 
which is, of course, the same as $z_i$ cabled by the skein element $v$
defined in Section~\ref{skein.sec1}.

Using the 
module
structure on $\BS(\Si)$ discussed above, we can now define the elements of our basis $\BB$ as follows:
\begin{eqnarray}\label{basedef} \bb(a,b,c) &=& v_1^{b_1} \cdots v_g^{b_g}
 \,\bb(a,0,c)\\
\nonumber &=& h^{-\sum_i b_i - \lfloor \frac 1 2 ( -e +\sum_i a_i ) \rfloor } (2+z_1)^{b_1}\cdots (2+z_g)^{b_g} \,\bg(a,0,c)~. \end{eqnarray}

\begin{thm} \label{basis}$\BB$ is a basis of $\BS(\Si)$.
\end{thm}

The proof will be completed in Section~\ref{proofsharp}. Note that it
is enough to show that these vectors lie in $\BS(\Si)$ and generate it
over $\BO$. As their number is equal to the dimension of 
$V_p(\Si)$, they must then form a basis, and $\BS(\Si)$ must be a free 
lattice. 

We call $\BB$ a {\em graph-like} basis. Note that $\bb(a,b,c)$ 
is a linear combination of vectors $\bg(a,b',c)$ with $b'_i\leq b_i$
for all $i=1,\ldots, g$. Moreover, the multiplicative expression in Eq.~(\ref{basedef}) shows that the
linear combinations are taken independently in each handle. Observe
also that the rescaling factor is a power of $h$ whose exponent
depends on the 
numbers 
$a_i$ and $b_i$ and the trunk color $2e$, but
does not explicitly depend on the colors $c_i$.

In genus one with no colored points, the basis $\BB$ is
 up to units
the same as the second integral basis $\{1,v,\ldots,v^{d-1}\}$ of \cite{GMW}. 
By functoriality, 
it follows 
that the vectors $v_i$ and their products lie in $\BS(\Si)$. However, it is not obvious {\em a priori} that the vectors $\bb(a,0,c)$ lie in $\BS(\Si)$. (Even in genus two the basis $\BB$  is different from the one obtained in \cite{GMW}.)

\begin{cor} The extended mapping class group 
(see Section~\ref{mcg.sec})
acts on $\BS(\Si)$ in the basis $\BB$ by matrices with coefficients 
in the cyclotomic ring $\BO$. Moreover, these matrices preserve a 
non-degenerate $\BO$-valued hermitian form.
\end{cor}

{\em The basis of $\BSs(\Si)$.}  
Recall that $\BSs(\Si)\subset V_p(\Si)$ is the dual lattice to $\BS(\Si)$ with respect to the $\BO$-valued hermitian form $(\ ,\ )_\Si$. (We will review the 
definition of this  form in Section~\ref{pB}.) Our basis of $\BSs(\Si)$ will be
denoted by $\BB^{\sharp}$. Using 
the 
algebra/module
structure  discussed above, we 
define the elements of  $\BB^\sharp$ as follows:
\begin{equation}\label{basedefsharp} \bb^{\sharp}(a,b,c) = h^{-\sum_i b_i - \lceil \frac 1 2 ( e+\sum_i a_i  ) \rceil } (2+z_1)^{b_1}\cdots (2+z_g)^{b_g} \,\bg(a,0,c)~. \end{equation}
as $(a,b,c)$ vary over the same index set. Note that the only
difference between Eqs.~(\ref{basedef}) and (\ref{basedefsharp}) is
that the trunk half-color $e$ appears with the opposite sign, and
$\lfloor \ \rfloor$ is replaced with $\lceil \  \rceil$. Thus 
$\bb^{\sharp}(a,b,c)$ is a rescaling of $\bb(a,b,c)$.

\begin{thm}\label{sharpbasis} $\BB^\sharp$ is a basis of $\BSs(\Si)$.
\end{thm}

The following sections until Section~\ref{proofsharp} will be devoted
to the proofs of Theorems~\ref{basis} and~\ref{sharpbasis}.

\section{\   The $3$-ball lemma.}\label{B3.sec}

Let $K_\BO(D^3, \ell_{n, 2})$ denote the skein module of $D^3$
relative $n$ points in $S^2=\partial D^3$ colored $2$. Here, the
notation $K_\BO$ means that we consider skein modules with
coefficients in $\BO$, not in $\BO[\frac 1 p]$. (But the result is
actually more general; see the remark at
the end of this section.) The following 
result will be needed in 
proving that the vectors $\bb(a,b,c)$ lie in 
$\BS(\Si)$.

\begin{thm}[$3$-ball lemma]\label{skeind3} If $n$ is even, $K_\BO(D^3,  \ell_{n, 2})$ is generated by skein elements which can be represented by a collection of $n/ 2$ disjoint arcs colored $2$. If $n$ is odd,  $K_\BO(D^3,  \ell_{n, 2})$ is generated by skein elements  which can be represented by a collection of $ (n-3)/ 2 $ disjoint arcs colored $2$ union one Y shaped component also colored $2$.
\end{thm}

Recall that the module $K_\BO(D^3, \ell_{n, 2})$ is generated by
    colored graphs which meet the boundary nicely in the given $n$
    points colored $2$. We remark here that the required Jones-Wenzl
    idempotents are all defined over $\BO$, because the only
    denominators needed in their definition 
are
the quantum integers
    $[n]$ ($n\leq p-2$) which
    are invertible in $\BO$.

The proof of Theorem~\ref{skeind3} proceeds by a series of lemmas.
The first is an exercise in skein theory whose proof is left to
 the reader (represent elements of $K_\BO(D^3, \ell_{n, 2})$ by diagrams in a disk and apply the usual fusion  formulas \cite{KL,MV}).

\begin{lem}\label{L0} $K_\BO(D^3, \ell_{n, 2})$ is generated by unions of tree graphs
where all edges are colored~$2$.
\end{lem}

In the remainder of this section,  unless otherwise stated,  unlabeled arcs in figures 
are assumed to be colored $2$. As usual, we put $\delta = -A^2-A^{-2}$. The
following two lemmas are simple skein-theoretical calculations.

\begin{lem}\label{L1} $ (A^4 -1 + A^{-4})
 \hskip.1in
 \begin{minipage}{0.2in}\includegraphics[width=0.2in]{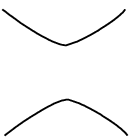}\end{minipage}
= A^{-4} \delta 
\hskip.1in \begin{minipage}{0.2in}\includegraphics[width=0.2in]{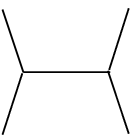}\end{minipage}
+ (1-A^{-4}) \ \delta
\hskip.1in \begin{minipage}{0.2in}\includegraphics[width=0.2in]{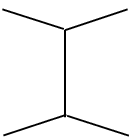}\end{minipage}
+
\hskip.1in \begin{minipage}{0.2in}\includegraphics[width=0.2in]{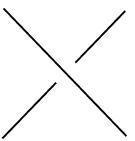}\end{minipage}.$
\end{lem}

\begin{proof}  We use the abbreviation:
  $\begin{minipage}{0.2in}\includegraphics[width=0.2in]{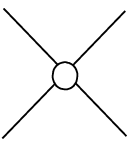}\end{minipage} = \begin{minipage}{0.2in}\includegraphics[width=0.2in]{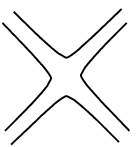}\end{minipage}$, where the arcs in the right hand diagram are colored $1$.
To prove \ref{L1}, it is enough to expand the right hand side using
\[\begin{minipage}{0.2in}\includegraphics[width=0.2in]{graphy/hibar.eps}\end{minipage}
  = \hskip.1in   
\begin{minipage}{0.2in}\includegraphics[width=0.2in]{graphy/x.eps}\end{minipage}
   - \delta^{-1}  
\begin{minipage}{0.2in}\includegraphics[width=0.2in]{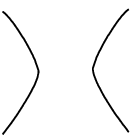}\end{minipage}
   \quad \text{and \cite[p.35]{KL}} \quad
\begin{minipage}{0.2in}\includegraphics[width=0.2in]{graphy/L+.eps}\end{minipage}
  = A^4 \hskip.1in   
\begin{minipage}{0.2in}\includegraphics[width=0.2in]{graphy/L0.eps}\end{minipage}
  + A^{-4} \hskip.1in   
\begin{minipage}{0.2in}\includegraphics[width=0.2in]{graphy/Linf.eps}\end{minipage}
  -\delta \hskip.1in   
\begin{minipage}{0.2in}\includegraphics[width=0.2in]{graphy/x.eps}\end{minipage}
  \]
\end{proof}

\begin{lem}\label{L2} $ 
 \begin{minipage}{0.2in}\includegraphics[width=0.2in]{graphy/vibar.eps}\end{minipage}
= (A^{4}-1) \delta^{-1} 
\hskip.1in \begin{minipage}{0.2in}\includegraphics[width=0.2in]{graphy/L0.eps}\end{minipage}
+ A^{-4} \delta^{-1}
\hskip.1in \begin{minipage}{0.2in}\includegraphics[width=0.2in]{graphy/Linf.eps}\end{minipage}
- \delta^{-1}
\hskip.1in \begin{minipage}{0.2in}\includegraphics[width=0.2in]{graphy/L+.eps}\end{minipage}.$
\end{lem}

\begin{proof} 
Rewrite the right hand side  and simplify using the equations in the proof of Lemma \ref{L1}:
\[ \delta^{-1} \left( A^4 \hskip.1in \begin{minipage}{0.2in}\includegraphics[width=0.2in]{graphy/L0.eps}\end{minipage} + A^{-4} \hskip.1in \begin{minipage}{0.2in}\includegraphics[width=0.2in]{graphy/Linf.eps}\end{minipage} - \begin{minipage}{0.2in}\includegraphics[width=0.2in]{graphy/L+.eps}\end{minipage} \right) - \delta^{-1} \hskip.1in \begin{minipage}{0.2in}\includegraphics[width=0.2in]{graphy/L0.eps}\end{minipage} =   
\begin{minipage}{0.2in}\includegraphics[width=0.2in]{graphy/x.eps}\end{minipage} -   \delta^{-1} \hskip.1in \begin{minipage}{0.2in}\includegraphics[width=0.2in]{graphy/L0.eps}\end{minipage} =   \begin{minipage}{0.2in}\includegraphics[width=0.2in]{graphy/vibar.eps}\end{minipage}.\]
\end{proof}

We refer to the diagram on the left hand side of Lemma \ref{L2}  as an
I-bar. Using this lemma repeatedly to expand I-bars in tree graphs colored $2$,
we see from Lemma~\ref{L0} that 
$K_\BO(D^3,  \ell_{n, 2})$ is generated by skein elements  which can
be represented as disjoint unions of arcs colored $2$ and $Y$-shaped
graphs colored $2$. Here we use that $\delta=-[2]$ is invertible in $\BO$. 

The crucial step is now the following lemma which
shows how to replace two $Y$-shaped
graphs colored $2$ by three arcs; clearly this will be enough to
complete the proof of Theorem~\ref{skeind3}. 
\begin{lem} The element \ \   
$ 
 \hskip.1in
 \begin{minipage}{0.4in}\includegraphics[width=0.3in]{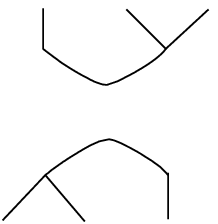}\end{minipage}$
 of $K_\BO(D^3, \ell_{6,2})$ is an $\BO$-linear combination of diagrams
 consisting of three arcs colored $2$. \end{lem}

\begin{proof}
Extending  all the diagrams in Lemma \ref{L1} by the same wiring, we obtain:
\[ (A^4 -1 + A^{-4})
 \hskip.1in
 \begin{minipage}{0.3in}\includegraphics[width=0.3in]{graphy/L0ex.eps}\end{minipage}
= A^{-4} \delta 
\hskip.1in \begin{minipage}{0.3in}\includegraphics[width=0.3in]{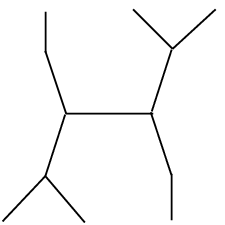}\end{minipage}
+ (1-A^{-4}) \delta
\hskip.1in \begin{minipage}{0.3in}\includegraphics[width=0.3in]{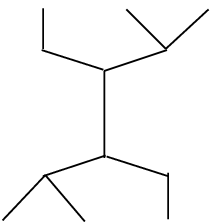}\end{minipage}
+
\hskip.1in \begin{minipage}{0.3in}\includegraphics[width=0.3in]{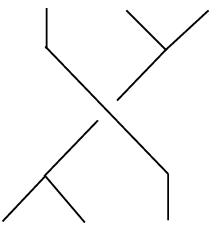}\end{minipage}.\]

Using Lemma \ref{L2} to expand two I-bars in the first diagram on the right hand side, two I-bars in the second diagram on the right hand side,
and one I-bar in the third diagram on the right hand side, we see that
the right hand side is an $\BO$-linear combination of diagrams
 consisting of $3$ arcs colored $2$. But $A^4 -1 + A^{-4}$ (which, up
 to units, is the sixth cyclotomic polynomial in $A^4$) is easily
 seen to be invertible 
in $\BO$. This proves the Lemma, and completes
 the proof of Theorem~\ref{skeind3}. 
\end{proof}

\begin{rem}{\em We never used that $A$ is a root of unity in this
    proof. Thus the result also holds for the skein
    module with coefficients in $\BQ(A)$, the ring of rational
    functions in $A$. In fact, it would be enough to work with the
    subring $\BZ[A,A^{-1}]$ with inverses of the relevant quantum
    integers and of $A^4 - 1 +A^{-4}$ adjoined to it. 
}\end{rem}

\section{\ Proof that $\bb(a,b,c)\in \BS(\Si)$.}

\begin{figure}[h]
\includegraphics[width=2in]{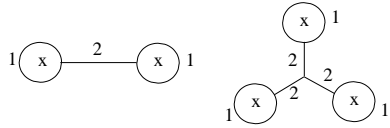}
\caption{An eyeglass and a tripod} \label{eyetri}
\end{figure}

 We fix some notation and terminology. Let $\Si_g$ denote the boundary of the handlebody $H_g$ with no 
colored points. 
The colored graph  that represents 
$\bg((1,1),(0,0))
\in \BS(\Si_2)$ is called an {\em eyeglass} and the colored graph  that
represents 
$\bg((1,1,1), (0,0,0))\in \BS(\Si_3)$ is called a
{\em tripod.} See Figure~\ref{eyetri}.  

\begin{lem}\label{eyeglasses} Eyeglasses and tripods are  divisible by $h$ in $\BS$.
\end{lem}

Since by definition $\bb((1,1),  (0,0))=h^{-1} \bg((1,1), (0,0))$ and $\bb((1,1,1), (0,0,0))=h^{-1} \bg((1,1,1),
(0,0,0))$, it follows that  $\bb((1,1), (0,0))
\in \BS(\Si_2)$ and 
$\bb((1,1,1), (0,0,0))
\in \BS(\Si_3)$.

\begin{proof} Recall $z_{i}$ denotes a simple loop which encloses the
$i$th hole. Let  $z_{i,j}$ denote a simple loop which encloses just
the $i$th and $j$th hole, etc. If we want these curves colored 
$v$, we just change the $z$ to a  $v.$ A scalar denotes this
scalar times the empty 
link.

If we expand the 
idempotent
$f_2$ in $\bg((1,1),(0,0))$,  we get $z_{1,2} +
[2]^{-1}
z_1z _2
.$  Making the substitution  $z= hv-2$, we get 
$$\bg((1,1),(0,0))= h v_{1,2} +  [2]^{-1} \left(h^2
v_1v _2 -2 h v_1 -2h  v_2 -2 \zeta_p^{-1} h^2  \right)~.$$ 
(We have used $2-[2]=2-\zeta_p-\zeta_p^{-1}=-\zeta_p^{-1}h^2$.)
As $[2]$ is a unit  of $\BO$ and $v$-graphs are in 
$\BS$
(Proposition \ref{vgraph}), this shows $\bg((1,1),(0,0))$ is 
divisible by $h$ in 
$\BS(\Si_2)$.   The divisibility of tripods is proved in the same way.
\end{proof}

Now consider again an arbitrary connected surface $\Si$, possibly with colored points $\ell(\Si)$.
Fix a lollipop tree $G$ and consider the elements $\bb(a,b,c)$ defined 
in Section~\ref{sec.thm}.

\begin{prop}\label{62} One has $\bb(a,b,c) \in \BS(\Si)$.
\end{prop}
\begin{proof} As already observed in Section~\ref{sec.thm}, it is enough
 to show this for $b=0$, since  $\bb(a,b,c)$ is obtained from $\bb(a,0,c)$ 
by multiplying it by some $v$-colored curves. In other words, we need to 
show that $\bg(a,0,c)$ is divisible by $h^{\lfloor \frac 1 2 ( -e +\sum_i a_i) \rfloor }$ in $\BS(\Si)$.

The colored graph representing $\bg(a,0,c)$ contains $g$ lollipops, where we
 mean by lollipop a subgraph consisting of a colored loop 
 joined
  to a stick at a single point. We now perform
a local change at the $g$ lollipops and obtain a new skein element $w(a,0,c),$
as is done in Figure \ref{w}. To do this, we can view the $i$th lollipop in $\bg(a,0,c)$ as $a_i$ arcs starting and ending at the 
idempotent $f_{2a_i}$ and looping around the $i$th hole. Here we refer to the usual device of representing an idempotent by a rectangle. Specifically, the $j$th arc connects the $j$th and the 
 $(2a_i-j+1)$th
input of the idempotent. The modification from $\bg(a,0,c)$ to $w(a,0,c)$
consists of inserting a braid so that  the  $a_i$ arcs now connect consecutive inputs on the $f_{2a_i}$. By a well-known property of the idempotent $f_{2a_i}$, this changes a given skein element only 
by multiplying it
by some power of $A$.
Hence it is enough to show that $w(a,0,c)$ is divisible by
the above-mentioned power of  $h$.
 \begin{figure}[h]
\includegraphics[width=3in]{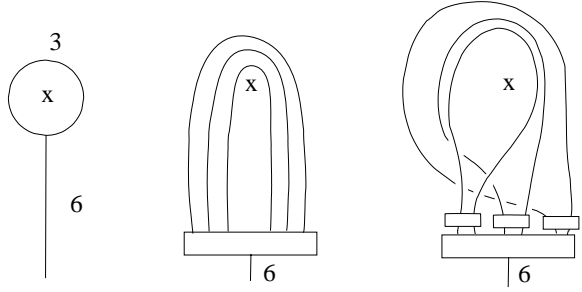}
\caption{a lollipop in  $\bg(a,0,c)$, its expansion as a skein diagram, and the corresponding portion of $w(a,0,c)$} \label{w}
\end{figure}

We can draw $w(a,0,c)$ slightly differently: we also insert $a_i$ idempotents $f_2$ (as is also done in Figure~\ref{w}).  This does not change the skein element at all (the $f_2$'s are ``absorbed'' by the $f_{2a_i}$ by another well-known property of the idempotent). We refer to this last operation as {\em spawning off} $f_2$'s from $f_{2a_i}$. If $\ell(\Si)$ is nonempty, we also spawn
 off $e$ 
idempotents $f_2$ 
  from the $f_{2e}$ on the trunk in the diagram for  $w(a,0,c)$.

There is a 3-ball  in $H_g$ whose boundary intersects $w(a,0,c)$ in
$e+ \sum_i a_i  $ points colored $2$ corresponding to the idempotents we spawned off. 
By the $3$-ball lemma~\ref{skeind3}, 
we can replace this part of the diagram with a linear combination over $\BO$ of diagrams each with
$\lfloor  \frac 1 2 ( e+ \sum_i a_i  ) \rfloor $ arcs  and Y's colored two.   At most $e$ of 
these arcs and Y's meet the trunk. Thus the rest of the arcs and Y's are completed to eyeglasses and tripods in the larger diagram. Thus $w(a,0,c)$ is also represented as a  linear combination over $\BO$ of diagrams
 with
  \begin{center}
$ \lfloor  \frac 1 2 ( e+ \sum_i a_i  ) \rfloor - e= 
   \lfloor  \frac 1 2 (-e + \sum_i a_i ) \rfloor  $
\end{center}
 eyeglasses and tripods. 
As each eyeglass and tripod is divisible by $h$ this gives exactly the required divisibility of $w(a,0,c)$. This completes the proof.
\end{proof}

\section{\   The lollipop lemma.}

The lollipop lemma (Theorem~\ref{ll} below) 
will be used to show that the elements $\bb^\sharp(a,b,c)$ lie in  
$\BSs(\Si)$.
Recall that a $v$-graph in
a $3$-manifold is a usual colored graph together with some banded 
link components colored $v$, where $v=h^{-1} (2+z)$ 
(see Section~\ref{skein.sec1}). 
As before, a subgraph  in a $v$-graph 
consisting of a colored loop  meeting an edge (called the stick)
is called a lollipop. A stick can take part in two lollipops. The color 
of the loop edge is at least one half of the 
stick color, but is allowed to be greater than that. Note that we have imposed no condition on how the  loop edge of a lollipop is embedded into the ambient manifold.

\begin{thm} \label{ll}(Lollipop Lemma) Let $L$ be a $v$-graph in $S^3$
  containing $N$ 
lollipops with stick colors $2a_1, 2a_2 \cdots 2a_N,$ 
then its evaluation $\langle L \rangle$ is divisible by $h^{\lceil \frac 1 2 {\sum_{i=1}^N a_i} \rceil }.$
\end{thm}

Here, the evaluation  $\langle L\rangle $ of a skein element $L$ in
$S^3$ is defined to be the ordinary Kauffman bracket.
Thus the empty link evaluates to $1$ and a
zero-framed 
unknotted loop colored one
evaluates to $ -\zeta_p -\zeta_p^{-1}$, for example. As shown in \cite{GMW}, the evaluation
of a $v$-graph in $S^3$ lies in $\BO$. 

Let us first prove a special case. A {\em basic lollipop} is a lollipop where the stick is colored $2$ and the loop edge is colored $1$.

\begin{lem}\label{6.1} The evaluation of a $v$-graph in $S^3$ with $N$ basic
  lollipops is divisible  by $h^{\lceil  N/ 2 \rceil}$ in $\BO$.
\end{lem}

\begin{proof} If we have a basic lollipop whose loop spans a disk which misses the rest of  the $v$-graph $L$ then
$\langle L\rangle =0$, by a well-known property of the Jones-Wenzl 
idempotents. More generally, if  $L$  intersects a 3-ball only in a basic lollipop, then $\langle L\rangle =0$.
Let us show that we may reduce to this case by changing crossings using the relation 
\[ 
\begin{minipage}{0.2in}\includegraphics[width=0.2in]{graphy/L+.eps}\end{minipage}
-  
\begin{minipage}{0.2in}\includegraphics[width=0.2in]{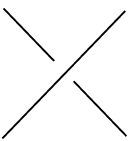}\end{minipage}
=   (A -A^{-1})
\left( \hskip.1in 
\begin{minipage}{0.2in}\includegraphics[width=0.2in]{graphy/L0} \end{minipage}
- 
\begin{minipage}{0.2in} \includegraphics[width=0.2in]{graphy/Linf.eps}  \end{minipage}
 \hskip.1in  \right) 
 \]
 and all the error terms are divisible by  $h^{\lceil N/ 2
  \rceil}$. 
(In this figure, strands are ordinary strands, {\em i.e.} colored
$1$.) 

By expanding all idempotents except one at the stick of each lollipop, we may assume that all arcs  of the graph are colored $1$.  When we change a crossing  between a $v$-colored link component and a strand colored $1$,  
the skein relation shows that the error term 
is given by 
$h^{-1}(A- A^{-1} )  =-A^{-1}$
times the difference  of two evalutions of $v$-graphs, each of which
satisfies again the hypothesis of the lemma 
but 
has
one less $v$-colored link component. By induction it follows  that  we can assume that  all $v$-colored link components are unlinked from the rest of the graph. Since a $v$-colored link evaluates to something integral, this shows that we may assume that there are no $v$-colored link components.

Next if we change crossings between a loop colored $1$ of a basic lollipop and any arc colored $1$, the error term is 
$-A^{-1}h$
times a $v$-graph with at least $N-2$ basic lollipops. By induction on the number of lollipops the error term satisfies  the conclusion of the lemma. Hence we may reduce to the case where one basic lollipop is not  linked with anything.
\end{proof}
\begin{rem}\label{altv}{\em In a similar way, we can give a new proof 
of the fact that the evaluation
of a $v$-graph lies in $\BO$. Namely, one checks that it is true for $v$-colored unlinks and then reduces to this case 
by
changing crossings, observing that error terms always lie in $\BO$.
}\end{rem}

\begin{proof}[Proof of Theorem \ref{ll}] We may expand the loop idempotents over $\BO$ into terms where each stick idempotent has $a_i$ arcs that meet it at $2a_i$ points along an edge of the idempotent.  Here we refer again to the device of representing an idempotent by a rectangle. As in the previous section, 
we insert  
a braid
in the strands that meet this edge so that the arcs now join points to their immediate neighbors.  This only changes the evaluation by a power of  $A.$ 
Then without changing the evaluation
one can spawn  from each $2a_i$-stranded idempotent  $a_i$ 
2-stranded idempotents.  Now  each term is a $v$-graph with  $\sum_{i=1}^N a_i$ basic lollipops
so the result follows from Lemma~\ref{6.1}.
\end{proof}

\section{\   Proof that  $\bb^\sharp(a,b,c)\in \BSs(\Si)$.}\label{pB}
In this section, we use the Lollipop Lemma to show that the elements 
 $\bb^\sharp(a,b,c)$ defined
in Section~\ref{sec.thm} lie in $\BSs(\Si)$. Let us first review
the definition of this module in more detail.

The hermitian form $( \ ,\ )_{\Si}$ is simply a rescaling of the
canonical hermitian form $\langle \ ,\ \rangle_{\Si}$ on $V_p(\Si)$
(see Section~\ref{skein.sec1}):
$$(x,y )_{\Si} =  \BD^{\beta_0(\Si)} \langle x,y\rangle_{\Si}~.$$
Here 
$\beta_0(\Si)$ is the number of components of $\Si$, and  $\BD$ is 
defined in Eq.~(\ref{defD}) in Section~\ref{sec.intro}. In fact, $\BD$ is the inverse of the quantum invariant $\langle S^3
 \rangle_p$.

When restricted to $\BS(\Si) \subset V_p(\Si)$, the form $(
 \ ,\ )_{\Si}$ takes values in $\BO$. This follows from the definition
 of $\BS(\Si)$ and the integrality result for quantum invariants of
 closed 
connected manifolds
\cite{Mu2,MR}. (See \cite{G,GMW} for more details.)

\begin{de} {\em The lattice $\BSs(\Si)$ is the dual lattice to
 $\BS(\Si)$ with respect to the hermitian form
$( \ ,\ )_{\Si}$.
}\end{de}

Now assume $\Si$ is connected. If $M$ has  boundary $\Si$, and $M$ is
also connected, 
let $K_v(M, \ell(\Si))$ denote the $\BO$-submodule of $K(M, \ell(\Si))$
spanned by $v$-graphs in $M.$ 
Proposition~\ref{vgraph} is equivalent to saying that for connected
$M$ and $\Si$, the natural map $$K_v(M, \ell(\Si)) \rightarrow \BS(\Si)$$
is surjective. 

\begin{rem}{\em The inclusions $K_\BO(M, \ell(\Si)) \subset K_v(M, \ell(\Si))\subset K(M, \ell(\Si))$ are strict in general.
}\end{rem}

The following Proposition~\ref{testingSsharp} is a useful device for
describing $\BSs(\Si)$. Let $H$ and $H'$ be two complementary handlebodies in $S^3$ with 
$\partial H = \Si = -\partial H'$. We define  a bilinear
form

\begin{center}
$( ( \ ,\ ) )_{H, H'}: K(H,\ell(\Si)) \times K(H',-\ell(\Si)) 
 \ \longrightarrow \ \BO[\frac 1 h]$
\end{center}   
by
\[  ((L, L'))_{H, H'} =   
\langle L \cup_{\ell(\Si)} L'\rangle. \]

Here $\langle \ \rangle$ is the usual Kauffman bracket of a colored
 graph in $S^3$. When restricted to 
$K_v(H,\ell(\Si)) \times K_v(H',-\ell(\Si))$, this form takes values in $\BO$.

 \begin{prop} \label{testingSsharp}
A skein element 
  $x \in K(H,\ell(\Si))$ represents an element of $\BSs(\Si)$ if and only if $((x,x'))_{H,H'} \in \BO$ for all $x' \in K_v(H',
-\ell(\Si))$. 
\end{prop}
\begin{proof}  The proof is essentially a standard argument in the skein-theoretical
  approach to TQFT's. 
The skein element $x$ defines the vector $[(H,x)]\in V_p(\Si)$.
The hermitian form 
is given by 
\begin{equation}\label{revers}
([(H,x)],[(H,y)])_\Si=
\BD \langle (H\cup_\Si -H, x \cup_{\ell(\Si)} 
y^\star
)\rangle_p~,
\end{equation}
where $y^\star$ denotes the skein element in $-H$ obtained from $y$ by reversing orientation.
This can also be viewed as a bilinear pairing of 
$[(H,x)] \in V_p(\Si)$ with $[(-H,y^\star)] \in V_p(-\Si)$.
But  
$[(-H,y^\star)]$ 
can also be represented by some skein element $y'\in K(H',-\ell(\Si))$. Thus 
 the hermitian pairing (\ref{revers}) is equal to $$
\BD \langle (H\cup_\Si H', x \cup_{\ell(\Si)} y')\rangle_p
=\langle x \cup_{\ell(\Si)} y'\rangle
=((x,y'))_{H,H'}$$ since $H\cup_\Si H'=S^3$. Now 
$$[(H,y)] \in \BS(\Si)\   \iff \ \
[(-H,y^\star)]
\in \BS(-\Si) \  
\iff 
\ y'\in   
K_v(H',-\ell(\Si))~.$$
Thus one has  
$[(H,x)]
\in \BSs(\Si)$ if and only if $((x,y'))_{H,H'} \in \BO$ for all $y' \in K_v(H',-\ell(\Si))$. 
\end{proof}

We are now ready to prove the following.

\begin{prop}\label{8.4} The elements $\bb^\sharp (a, b,c)$ lie in 
$\BSs(\Si).$
\end{prop}

\begin{proof}  Embed the handlebody $H_g$ into $S^3$ so that its exterior is also a handlebody 
 $H_g'$.
By Proposition \ref{testingSsharp}, we only need show that 
if $\bb^\sharp (a,b,c)$ in $H_g$ is completed by any $v$-graph in
$H_g'$,
then the evaluation of the result lies in $\BO$. However we
may isotope the 
$v$-colored curves
$v_i^{b_i}$ ($i=1,\ldots, g$) 
  in 
$\bb^\sharp (a,b,c)$ 
across $\Si.$ Thus we only need to prove this statement in the case $b=0.$ In other words, we only need to show that if $\bg(a,0,c)$ is completed by any $v$-graph in the complementary handlebody then the evaluation of the result is 
divisible by  $h^{ \lceil \frac 1 2 ( e+\sum_i a_i ) \rceil }$ in $\BO.$ 

If $e=0$, this follows immediately from the Lollipop Lemma~\ref{ll}, as in
any completion of $\bg(a,0,c)$ we have the $g$ lollipops of $\bg(a,0,c)$ 
with the sum of the stick half-colors equal to $\sum_i a_i.$

If $e>0$, then in any completion of $\bg(a,0,c)$ we first modify 
the part 
in $H_g'$ which is glued to 
the trunk edge 
of $\bg(a,0,c)$ 
in a now familiar way.
To do this, we represent the idempotent $f_{2e}$ on the trunk edge diagrammatically by the usual rectangle
and expand all the idempotents 
in the glued-on part of the graph 
into strands
 colored $1$. In every term of this expansion, we have $e$ arcs that start and end at the 
``bottom'' of this 
rectangle.
Inserting an appropriate braid we can arrange that the arcs now join consecutive points on the ``bottom'' of 
the rectangle, and then we spawn off
$e$ idempotents $f_2$. As before, we  can compensate for these changes by changing the  coefficients in the expansion over $\BO$ by some powers of $A$. In each term, we now see $e$ basic lollipops with stick color $2$
below the trunk edge. 
Using these
and the $g$ lollipops in $\bg(a,0,c)$ itself, we see that any completion of
$\bg(a,0,c)$ can be expanded over $\BO$ as a 
linear combination 
of $v$-graphs containing 
lollipops with the sum of the stick 
half-colors equal to
$e+\sum_i a_i.$ 
 By the Lollipop Lemma~\ref{ll} we are done.
\end{proof}

\section{\  Index counting 
}\label{proofsharp}

Let $\bB$ be the 
$\BO$-lattice in $V_p(\Si)$ spanned by the $\bb(a,b,c).$   Let $\bB'$ be the 
$\BO$-lattice 
spanned by the $\bb^\sharp(a,b,c)$. \footnote{We don't use the notation $\bB^\sharp$ because it is not clear {\em a priori} that $\bB'$ is the dual lattice of $\bB$.} We know $\bB \subset \BS(\Si)$
and $\bB' \subset \BSs(\Si)$ by 
Propositions~\ref{62} and~\ref{8.4}.
In this section, we will use order ideals to show that both inclusions are equalities (Theorem~\ref{9.3}).  
This will prove Theorems~\ref{basis} and~\ref{sharpbasis}. We will sometimes abbreviate  
$\BS(\Si)$ by $\Bs$ and  $\BSs(\Si)$ by $\Bss$.

Let $\bG$ denote the 
$\BO$-lattice
spanned by the small graph basis $\BG$ of $V_p(\Si)$, {\em i.e.}  
the elements $\bg(a,b,c)$. 
One has $\bG \subset \bB\subset \bB'$ and all three of these lattices are free. 

If $L \subset  L'$ is an inclusion of free lattices of the same rank over  $\BO$, define their {\em index} $[L':L]$ to be the determinant up to units in $\BO$ of a matrix representing a basis for $L$ in terms of a basis for $L'.$

\begin{prop} \label{numer}$[\bB:\bG] [\bB' :\bG] =  [\bG^\sharp : \bG]$~.
\end{prop}

\begin{proof}
By \cite[5.4]{GMW}, we have 
\[  [\bG^\sharp: \bG]= h ^{g(d-1)\dim(V_p(\Si))}. \] 
Actually \cite[5.4]{GMW} only deals with the case $\ell(\Si)= \emptyset,$ but this equation holds in general by the same argument.
By construction: $[\bB : \bG]=h^{N}$ where
$N= \sum_{(a,b,c)}\left( {\sum_i b_i + \lfloor \frac 1 2 ( -e +\sum_i a_i   ) \rfloor}\right)$.
Also by 
 construction: $[\bB' : \bG]=h^{N'}$ where
$N'= \sum_{(a,b,c)} \left( {\sum_i b_i + \lceil \frac 1 2 ( e+ \sum_i a_i  ) \rceil}\right).$ Of course, $e$ depends on $(a,b,c)$ in these expressions. Luckily, $e$ drops out when we look at $N+N'$. 

We need to show $N+N'= g(d-1)\dim(V_p(\Si)),$ or 
\[ \sum_{(a,b,c)} \left( 2\sum_i b_i +  \sum_i a_i   \right) = g(d-1)\sum_{(a,b,c)} 1. \]
It suffices to prove this with $a$ and $c$ held constant:
\[ \sum_{b}\ \left( 2\sum_{i=1}^g b_i +   \sum_{i=1}^g a_i   \right) = g(d-1)\sum_{b} 1. \] 
Recall that the set of $b$ being summed over is the set of $(b_1, \ldots, b_g)$
such that $0 \le b_i \le d-1-a_i$, and the cardinality of this set is 
$\prod_{j=1}^g (d-a_j)$.
So it suffices to show: \[ 2\sum_{b}\  \sum_{i=1}^g b_i +  \left( \sum_{i=1}^g a_i\right) \prod_{j=1}^g (d-a_j) =
g (d-1) \prod_{j=1}^g (d-a_j)  \] 
Noting that $g (d-1) = \sum_{i=1}^g ( d-1),$ we only need to show
\[ 2\sum_{b}\  \sum_{i=1}^g b_i =  \sum_{i=1}^g (d-1-a_i) \prod_{j=1}^g (d-a_j)  \] 
This equation expresses the fact that the average value of $\sum_{i=1}^g b_i$
over the index of summation $b$ is $ \frac 1 2  \sum _i (d-1-a_i). $ To see this, we consider the involution on the index set $\{b\}$ defined by sending $b_i$ to $d-1 -a_i -b_i$ for each $i.$ 
On each orbit the average value of $\sum_{i=1}^g b_i$
 is $ \frac 1 2  \sum _i (d-1-a_i).$ This establishes the identity.
  \end{proof}

If $M$ is a finitely generated torsion module over $\BO$, we denote its order ideal by 
$\oi(M).$  If $L \subset  L'$ is an inclusion of free lattices of the same rank over  $\BO$, then $\oi(L'/L)$ is the
principal
ideal 
in $\BO$
generated by 
$[L':L]$.

We will need the following result (here $\overline{M}$ denotes the conjugate module of $M$).

\begin{prop} \label{S/G} $\bG^\sharp / \Bss \cong \overline{\Bs/\bG}$~.
\end{prop}

Proposition~\ref{S/G} is a general fact about inclusions of 
(not necessarily free) lattices of the same rank 
equipped with a non-singular hermitian form over a Dedekind domain.  We were unable to find this fact stated in the literature, but it is not hard to deduce it from 
\cite[Theorem 4.12]{R}. We omit the details.

The proof of Theorems~\ref{basis} and~\ref{sharpbasis} is completed with the following 
\begin{thm}\label{9.3}  $ \BS(\Si)= \bB$ and  $ \BSs(\Si) =\bB'$.  In particular $\BS(\Si)$ and $\BSs(\Si)$ are free.
\end{thm}

\begin{proof}
As $\bG \subset \bB \subset \Bs$, we have that
\begin{equation} \label{E1} 
\oi(\Bs/\bG) = \oi(\Bs/\bB) \oi(\bB/\bG).
\end{equation}
Similarly $\bG \subset \bB' \subset \Bss$ gives us:
 \begin{equation} \label{E2}
 \oi(\Bss/\bG) = \oi(\Bss/\bB') \oi(\bB'/\bG).
 \end{equation}
 \begin{align*}
 \oi(\bG^\sharp / \bG) & = \oi(\bG^\sharp /  \Bss)  \oi(\Bss/\bG) \\  
                & = \overline{ \oi(\Bs /  \bG ) } \oi(\Bss/\bG) \text{ by Proposition \ref{S/G} } \\   
        & = \overline{ \oi(\Bs/\bB)} \overline{ \oi(\bB/\bG)}  \oi(\Bss/\bB') \oi(\bB'/\bG)
         \text{ by Equations (\ref{E1}) and    (\ref{E2})} \\   
        & = \overline{ \oi(\Bs/\bB)} \oi(\bB/\bG)  \oi(\Bss/\bB') \oi(\bB'/\bG) \text{ as $h$ is 
self-conjugate 
up to units} \\   
        & = \overline{ \oi(\Bs/\bB)}   \oi(\Bss/\bB')\oi(\bG^\sharp / \bG)  \text{ by Proposition \ref{numer}}.  
         \end{align*}
        Thus $ \oi(\Bs/\bB)= \oi(\Bss/\bB')= 1, $ hence $\Bs= \bB$ and $\Bss=\bB'.$
\end{proof}

\section{The necessity of a trunk}

In this section we give an example which explains why we have insisted that lollipop trees have trunks.  Let $\Si$ be the boundary of a regular neighborhood of a lollipop tree $G.$  A graph-like basis for $G$ is a basis for $\BS(\Si)$  that is obtained from the  small graph basis for $\Vp(\Si)$ by ``peeling'' off any excess $b$ from each loop color of a small basis vector, inserting a compensating $(z+2)^b$ and  then rescaling the resulting elements by some factor. 

\begin{figure}[h]
\includegraphics[width=2.1in]{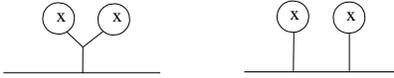}
\caption{ The graphs $G$ and $G'$} \label{gg}
\end{figure}

Now let $\Si$ denote a surface of genus two with two banded points colored $2$. Consider the graphs pictured in Figure \ref{gg} .  

While $G$ is a lollipop tree, $G'$ is not, as it does not have a trunk. Still, admissible colorings with small colors  on $G'$ also yield a basis for $V_p(\Si)$, and the notion of graph-like basis makes sense for $G'$. However, we have the following 

\begin{thm} $\Bs_5(\Si)$ does not have a graph-like basis associated to the graph $G'$.
\end{thm}

\begin{proof} 
Figure \ref{gs}  illustrates  two  basis elements from $G$ and two from $G'.$ In this figure,  the loops are all colored $1$ and the edges all colored $2$. 

 \begin{figure}[h]
\includegraphics[width=4in]{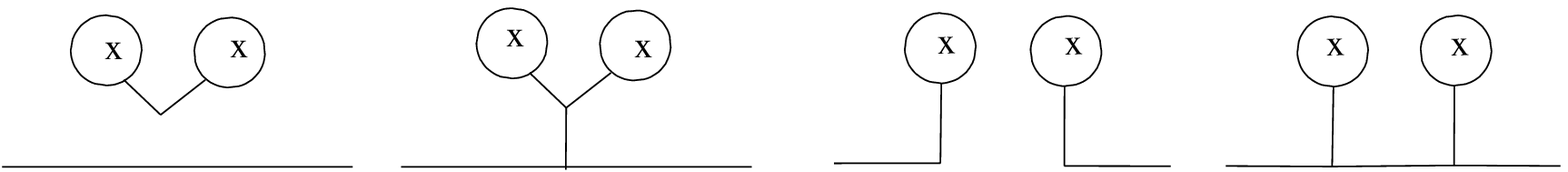}
\caption{ The basis elements $\bg_0$,   $\bg_1$, $\bg_0^\prime$ and $\bg_1^\prime$} \label{gs}
\end{figure}

Using the fact that 
$\BS$
of a 2-sphere with four  points colored $2$ is two-dimensional, one can see that  
\begin{align} \bg_0  &=  a \bg_0^\prime+ b \bg_1^\prime  \label{e4} \\
\bg_1  &=  \alpha \bg_0^\prime+ \beta \bg_1^\prime,  \label{e5}
\end{align} 
for some $a,b,\alpha, \beta \in \BO$. The exact values of these are not important, but we will need to use that $a,b\in \BO^\star$ which is easily checked.

 As $G$ is a lollipop tree, we know that \begin{equation}  h^{-1} \bg_0 \in \Bs_5(\Si) \text{ and }   h^{-1} \bg_1 \notin \Bs_5(\Si)~.  \label{e6}
\end{equation}

 Let us try to express $h^{-1} \bg_0$  in terms of a graph-like basis for $G'$,
assuming such a basis exists. We can modify the elements $\bg_0^\prime$ and $\bg_1^\prime$ only by a rescaling, as no peeling off is possible for these elements. Hence the required expression for $h^{-1} \bg_0$ would be   
$$h^{-1} \bg_0 =  a h^{-1} \bg_0^\prime+ bh^{-1} \bg_1^\prime~.$$ Since $a$ and $b$ are units of $\BO$, both $h^{-1} \bg_0^\prime$ and $h^{-1} \bg_1^\prime$ would have to exist in $\Bs_5(\Si)$.
But then Eq.~(\ref{e5}) implies that $h^{-1}\bg_1$ would also exist in $\Bs_5(\Si)$, which contradicts the second half of Statement~(\ref{e6}). Thus $\Bs_5(\Si)$ does not have a graph-like basis for $G'$. 
\end{proof}

\section{Integral modular functor}

The collection of 
$\BO[\frac 1 p]$-modules
$V_p(\Si)$ associated to surfaces with
colored points form a modular functor (see {\em
  e.g.} \cite{T}). This means 
that the
$V_p(\Si)$ satisfy certain
axioms describing their
behaviour when 
surfaces are cut into pieces in various ways. These axioms reflect 
the semi-simplicity of the underlying modular category \cite{T}. 

The integral lattices $\BS(\Si)$
should be a guiding example towards the notion of an ``integral 
modular functor''. But their behaviour under such ``cut and paste''  
operations is considerably more
complicated. In some sense, the reason is that the integral theory is 
no longer semi-simple.

 We hope to develop these ideas elsewhere. Here, we limit ourselves 
to describing
one axiom in a particularly simple situation, where a 
rescaling of the usual modular functor axiom suffices.

Let $\gamma$ be a separating simple closed curve on the connected
surface $\Si$. The curve $\gamma$ cuts $\Si$ into two subsurfaces
$\Si'$ and $\Si''$. If $\Si$ has colored points, we assume $\gamma$ 
misses these points. 
 For simplicity, we
 also assume that the sum of the colors on $\Si'$ (hence also on $\Si''$) is even.
Let $\Si'_i$ and $\Si''_i$ be the subsurfaces
with their boundary capped off by a disk containing an additional $i$-colored
 banded point.  
Then \cite[1.14]{BHMV2} the obvious gluing map along the disks induces an isomorphism 
\begin{equation}\label{IF1}
\bigoplus_{i=0}^{d-1} V_p(\Si'_{2i}) \otimes V_p(\Si''_{2i})
\mapright\approx V_p(\Si)~.
\end{equation}

This induces an injective map
\begin{equation}\label{IF2}
\bigoplus_{i=0}^{d-1} \BS(\Si'_{2i}) \otimes \BS(\Si''_{2i})
 \longrightarrow
\BS(\Si)~,
\end{equation}
whose image is a free sublattice of $\BS(\Si)$. 

\begin{thm}\label{11.1} If at most one of $\Si'$ and $\Si''$ has colored points,
  then there exist bases $\BB'_{2i}$ of $\BS(\Si'_{2i})$, $\BB''_{2i}$ of
  $\BS(\Si''_{2i})$, and $\BB$ of $\BS(\Si)$, such that $\BB$ is a
  rescaling of the tensor product basis $\sqcup_i \BB_{2i}'\otimes \BB_{2i}''$. (Here,
  $\sqcup_i \BB_{2i}'\otimes \BB_{2i}''$ is viewed as a basis of $V_p(\Si)$ via the map \eqref{IF1}.) 

More precisely, there are functions $\varphi_i:\BB'_{2i} \times \BB_{2i}'' 
\rightarrow
\BZ_{\geq 0}$ such that 
\begin{equation} \label{IF3} 
\BB = \sqcup_i \{ h^{-\varphi_i(\bb',\bb'')} \bb'\otimes \bb'' \,\vert\, \bb'\in \BB_{2i}',
\bb'' \in \BB_{2i}'' \}
\end{equation}
\end{thm}

\begin{proof} Write $\Si$ as the boundary of a handlebody $H$ such
  that the curve $\gamma$ bounds a disk $D$ in $H$. This disk cuts
  $H$ in handlebodies $H'$ and $H''$ with boundary 
$\Si'\cup D$ and $\Si''\cup D$.
We can find a 
graph-like basis $\BB$ of
  $\BS(\Si)$  with respect to a lollipop tree $G$ in $H$ such that $G$
  meets the disk $D$ transversely in one edge, 
and 
such that, moreover, $G'=G\cap H'$ is a lollipop tree in $H'$, and
$G''=G\cap H''$ is a lollipop tree in $H''$. Here we use the
hypothesis that at most one of $\Si'$ and $\Si''$ has colored points. 
Taking for $\BB_{2i}'$ and $\BB_{2i}''$ the graph-like bases associated to $G'$
and $G''$, with color $2i$ on the edge meeting the disk $D$,  the result follows.
\end{proof}

It is easy to write down an explicit formula for the 
rescaling factors $\varphi_i :\BB_{2i}' \times \BB_{2i}'' \rightarrow
\BZ_{\geq 0}$ in this situation; this is left to the reader.

\begin{rem} {\em (i) If both $\Si'$ and $\Si''$ have colored points, then
    in general 
there is no basis of $\BS(\Si)$ which is just a
    rescaling
of a tensor product basis as in \eqref{IF3}. This follows from the
    example in the previous section.

(ii) If $\gamma$ is non-separating, the modular functor axiom also needs 
more modification than just a rescaling.}
\end{rem}

\section{Disconnected surfaces and the tensor product axiom
}\label{discon.sec}

Let $\Si$ and $\Si'$ be  two closed surfaces.
We have compatible natural maps $\BS(\Si) \otimes \BS(\Si') \rightarrow \BS(\Si \sqcup \Si')$ and $\Vp(\Si) \otimes \Vp(\Si') \rightarrow \Vp(\Si \sqcup \Si')$
specified by sending
$[(M,L)] \otimes [(M',L')]$ to $[(M \sqcup M',L \sqcup L')]$ where $(M,L)$ and $(M',L')$ are $3$-manifolds with colored links whose  boundary is $\Si$ and $\Si'$, respectively.

The map  $\Vp(\Si) \otimes \Vp(\Si') \rightarrow \Vp(\Si \sqcup \Si')$
is an isomorphism \cite{BHMV2}, and this property is called the tensor
product axiom.
We are interested in the extent that $\BS(\Si) \otimes
\BS(\Si') \rightarrow \BS(\Si \sqcup \Si')$ is also an isomorphism.
It follows from Corollary \ref{tpa} 
below
that this map is an isomorphism if
$\Si$ and $\Si'$ have no colored points. If there are colored points,
the image of this 
map
 may only be a sublattice of $\BS(\Si \sqcup
\Si')$. But we shall see that a basis of 
$\BS(\Si \sqcup \Si')$ can always be obtained as a rescaling of a 
tensor product basis. 

The following definition allows for a convenient expression of the
needed rescaling.
 Let $\BB$ be 
  the
   basis of $\BS(\Si)$,
for a connected surface $\Si$,
 associated to some
lollipop tree $G$ in a handlebody $H$. 
We define the {\em oddity} $\varepsilon(\bb)$ of a
basis element $\bb\in\BB$ as follows. \footnote{A different notion of  
parity will be defined in Section~\ref{torsion.sec}.} Let $2a_1, \ldots, 2a_k$ be the
colors of the stick edges, and let $2e$ be the trunk color of 
$\bb$. Let $A(\bb) =\sum_i a_i $, and
$e(\bb) =e$.  Define 
$\varepsilon(\bb)=1$ if $e(\bb)=d-1$ and $A(\bb) -
e(\bb)$ is odd. Otherwise, $\varepsilon(\bb)=0$.

\begin{thm}\label{DC} 
Let $\Si_1, \ldots, \Si_n$ be  connected surfaces, and
let $\BB_1, \ldots,
  \BB_n$ be graph-like bases of $\BS(\Si_1), \ldots, \BS(\Si_n)$. Then
  the set
\begin{equation} \label{DC1} 
\BB = \{ h^{-\lfloor \frac 1 2  \sum_i \varepsilon(\bb_i)  \rfloor }
\bb_1\otimes \cdots \otimes \bb_n \,\vert\, \bb_i\in \BB_i \}
\end{equation}
is a basis of 
      $\BS(\sqcup_i \Si_i)$.
\end{thm}

\begin{figure}[h]
\includegraphics[width=2.5in]{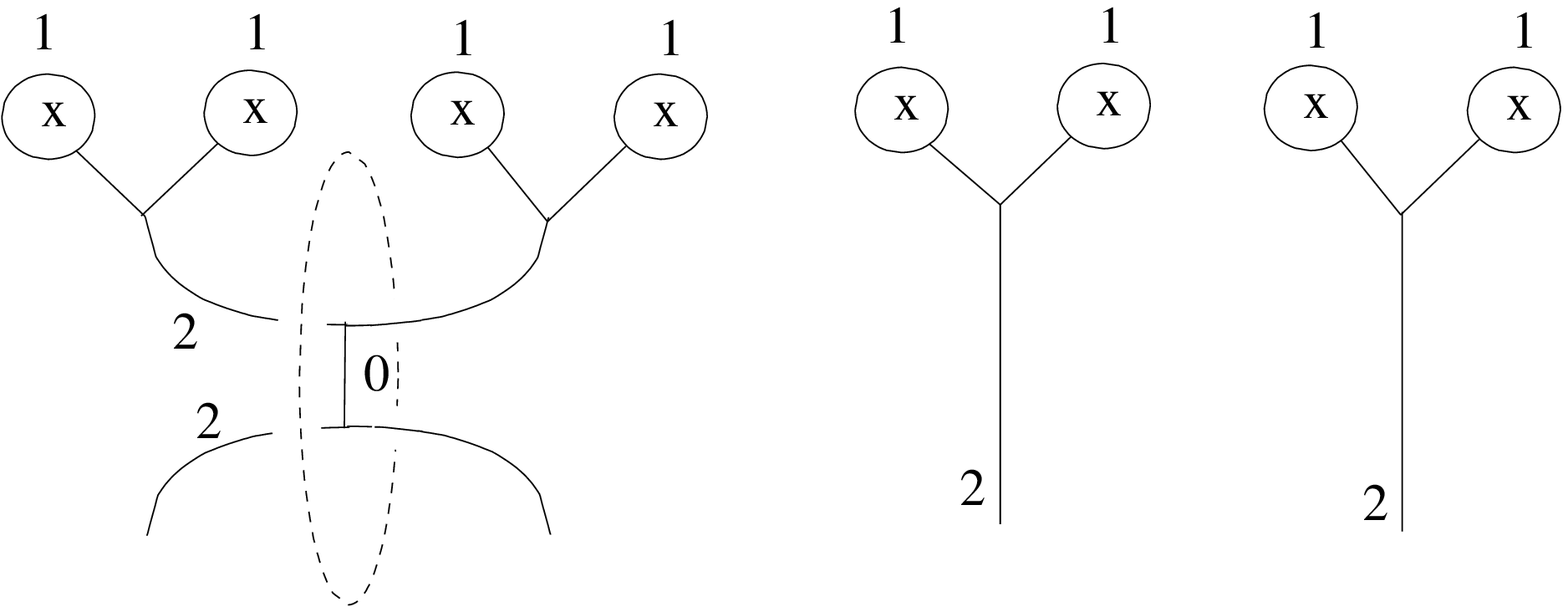}
\caption{} \label{exa}
\end{figure}

\begin{rem}{\em  If $\Si$ is connected, the lattice $\BS(\Si)$ has
    basis vectors with non-trivial oddity if and only if $\Si$ has
    genus at least two and the sum of the colors of the banded points
    on $\Si$ is at least $2(d-1)=p-3$.  For example, the rightmost diagram in 
    Figure~\ref{exa} represents an element with non-trivial oddity 
in $\Bs_5$ of a genus $2$ surface with one colored point colored $2$.}  
\end{rem}
 
\begin{cor} \label{tpa} The natural map $\BS(\Si) \otimes \BS(\Si') \rightarrow \BS(\Si \sqcup \Si')$
is always injective with a cokernel isomorphic to a direct sum of
cyclic modules $\BO/h\BO$. It is an isomorphism if one of the surfaces has no colored points.
  \end{cor}

\begin{rem}{\em In fact, the map $\BS(\Si)
    \otimes \BS(\Si') \rightarrow \BS(\Si \sqcup \Si')$ has a
    non-trivial cokernel if and
    only if both surfaces have a connected component whose
    $\BS$-lattice has basis vectors with non-trivial oddity. This
    follows  from Theorem \ref{DC}.}
\end{rem}

Since the tensor product axiom for $\BS$ does not always hold for surfaces with colored points, it seems that any proof of the tensor product axiom for surfaces without colored points must ultimately depend on detailed knowledge of bases for $\BS$ for connected surfaces.

For the proof of Theorem \ref{DC},  we need
to state the analog of Proposition \ref{vgraph} for disconnected
surfaces. 
The proof is the same as in the connected case.
A little terminology is convenient. Let $\pi_0(\Si)$ be the set of connected components of $\Si$. We say
that a $3$-manifold $M$ with boundary $\Si$ represents a partition $P$ of
$\pi_0(\Si)$ if $P$ coincides with the partition given by the fibers of the natural map
$\pi_0(\Si)\rightarrow \pi_0(M)$. 

\begin{prop} \label{vgraph2} If $\Si$ is not connected, choose any collection $\mathcal M$ of
  $3$-manifolds $M$ with boundary $\Si$ such that every partition of
$\pi_0(\Si)$ is represented by some $M\in \mathcal M$. Then $\BS(\Si)$ is
  generated over $\BO$ by elements in $V_p(\Si)$  represented by $v$-graphs in
  the $3$-manifolds $M\in \mathcal M$, where, as before, the 
  $v$-graphs 
  must meet the boundary in the colored points of
  $\Si$.
 \end{prop}

\begin{proof}[Proof of Theorem \ref{DC} in the case n=2]
 
We suppose $\Si_i$ is the boundary of $H_{g_i}.$  
We wish to find a basis for $ \BS(\Si_1 \sqcup \Si_2)$ which we recall is defined as a subset  of  $\Vp(\Si_1 \sqcup \Si_2)$ which we can identify  with $\Vp(\Si_1) \otimes \Vp(\Si_2)$. Thus any element of $ \BS(\Si_1 \sqcup \Si_2)$ may be described as a linear combinations over $\BO[\frac 1 h]$ of  elements of the form $\bb_1 \otimes \bb_2.$  Here $\bb_i$  is defined with respect to some lollipop tree $G_i$ in $H_{g_i}$.

Let $\#$ denote the operation of (interior)
connected sum of connected $3$-manifolds. $\Si _1\sqcup \Si_2$ is the boundary of $H_{g_1} \# H_{g_2}.$  According to Proposition \ref{vgraph2},  
$\BS(\Si_1 \sqcup \Si_2)$ is generated by $v$-graphs in $H_{g_1} \sqcup  H_{g_2}$ and  $v$-graphs in $H_{g_1} \# H_{g_2}$ .  The $v$-graphs in $H_{g_1} \sqcup  H_{g_2}$
generate the subspace of $\BS(\Si_1 \sqcup \Si_2)$ that has as basis the set of $\bb_1 \otimes \bb_2$ associated to small colorings of $G_1$ and $G_2$.  
We want to see if $v$-graphs in $H_{g_1} \# H_{g_2}$ can give any elements not in this subspace.

Let $\#_{\partial}$ denote the operation of boundary
connected sum of connected $3$-manifolds with connected boundaries.
$H_{g_1} \# H_{g_2}$ is homeomorphic to $H_{g_1} \#_\partial H_{g_2}$ with a 2-handle attached along the curve which bounds the disk used to form the boundary connected sum. We indicate such curves in figures as dotted curves. 
Thus dotted curves always indicate that a 2-handle should be attached to a handlebody along the curve.  The handlebodies are not actually drawn in our figures but are regular neighborhoods of the lollipop trees. We identify  $H_{g_1} \#_\partial H_{g_2}$ with  $H_{g_1 +g_2}$. 
\begin{figure}[h]
\includegraphics[width=2.5in]{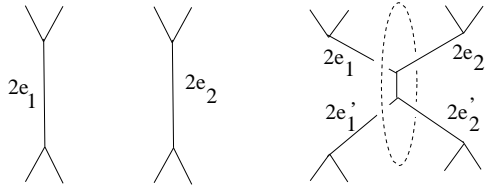}
\caption{ Lollipop trees $G_1$ in $H_{g_1}$, $G_2$ in $H_{g_2}$  and $G$ in $H_{g_1} \# H_{g_2}.$ We say $G$ is the grafting of $G_1$ and $G_2$  along their trunks (only the parts near the trunks are shown, and all loop edges are above the trunks in this figure).} \label{graft}
\end{figure}

A $v$-graph in $H_{g_1} \# H_{g_2}$ can be isotoped into $H_{g_1 +g_2}$ where we know a graph-like basis for $\BS(\partial (H_{g_1 +g_2}))$ associated to the lollipop tree $G$ obtained by the grafting of $G_1$ and $G_2$ along their trunks. See Figure \ref{graft}. 
 If
$\bb$ is such a basis element of $\BS(\partial (H_{g_1+g_2}))$,  we let $\hat \bb$ denote the element that $\bb$ represents in $\BS(\Si_1 \sqcup \Si_2 ).$ Our task is to compute $\hat \bb$.  

Suppose the coloring of $G$ that gives $\bb$ is as shown on the right of Figure \ref{graft}. Note that the 2-sphere composed of the disk spanning the dotted curve and the core of the 2-handle which is attached along this dotted curve can be isotoped to intersect the colored lollipop tree in either two points colored  $2e_1$ and  $2e'_1$ or two points colored  $2e_2 $ and $2e'_2$.  As $V_p$ of a 2-sphere with two points with distinct even colors is the zero module, we see that $\hat \bb$ is zero unless $e_1 = e'_1$ and $e_2 = e'_2$. In this  case, we may consider the small colorings of the $G_i$ which agree with the colorings of $G$ except on the edges already labelled $2e_i$ in Figure \ref{graft}. Let $\bb_i$  denote the basis elements of $\BS(\Si_i)$ indexed by these small colorings of $G_i.$ 
In this case,  
 $\hat \bb$ is up to units $ h^E\   \bb_1 \otimes  \bb_2$ where
\begin{center} $E= { d-1 -  \lfloor \frac 1 2 ( A(\bb)  -e(\bb)) \rfloor+ 
\lfloor \frac 1 2 ( A(\bb_1)  -e(\bb_1) ) \rfloor +  
\lfloor \frac 1 2 (  A(\bb_2)  -e(\bb_2)) \rfloor
}.$\end{center}
To see this one uses fusion and the fact that  $V_p$ of a $2$-sphere with a single point with a nonzero even color is the zero module. (When applying fusion, one encounters certain coefficients, but these are products of quantum integers and their inverses, hence units in $\BO$.) 
One also must take into account the powers of $h$ in the definitions of $\bb$,$\bb_1$, and $\bb_2.$ The $d-1$
term comes from the surgery axiom (S1) of  \cite{BHMV2} \footnote{In \cite{BHMV2}, $\BD$ is denoted by $\eta^{-1}$.} which implies
that attaching  a 1-handle to a $3$-manifold has the effect of multiplying its quantum invariant  by $\mathcal{D}$ which we recall, up to units, is $h^{d-1}.$

One checks that if E is negative, it must be  $-1$, and this  happens precisely when
$e(\bb)=0$, $e(\bb_1)=e(\bb_2)=d-1$, and $A(\bb_1) \equiv A(\bb_2) \equiv d  \pmod{2},$ in other words, when $\varepsilon(\bb_1)= \varepsilon(\bb_2)=1$. Conversely, if $\varepsilon(\bb_1)= \varepsilon(\bb_2)=1$ we can find a $\bb$ such that $\hat \bb$ is $h^{-1} \bb_1\otimes \bb_2$ up to units. 
Thus the elements $h^{-\lfloor \frac 1 2  (\varepsilon(\bb_1) + \varepsilon(\bb_2)) \rfloor }
\bb_1\otimes  \bb_2$ form indeed a basis of $\BS(\Si_1\sqcup \Si_2)$.
 \end{proof}

\begin{ex}\label{ex}{\em Let $\Si = \partial H_2$ with one point colored $2$. On the left of Figure \ref{exa} is pictured a skein element in $H_2 \# H_2$   representing an element of  $\Bs_5(\Si \sqcup \Si)$ which is divisible by $h^2$. On the right of the figure is a skein element in $H_2 \sqcup H_2$ representing an element of $\Bs_5(\Si \sqcup \Si)$ which (by the argument in the proof above) is $h$ times the element on the left (up to units) and thus is divisible by $h$ in $S_5(\Si \sqcup \Si)$. This latter element is of the form $\bb\otimes\bb$ where $\bb\in \Bs_5(\Si)$ is not divisible by $h$. This shows that the cokernel of the homomorphism from $S_5(\Si) \otimes  S_5(\Si) \rightarrow S_5(\Si \sqcup \Si)$ contains a non-trivial element which is annihilated by $h$. Of course, one has $\varepsilon(
\bb
)=1$.} \end{ex}

\begin{proof}[Proof of Theorem \ref{DC} for  $n > 2$]

As before, we suppose $\Si_i$ is the boundary of $H_{g_i}.$
Given a partition of $\pi_0( \sqcup_i \Si_i)$, we can take a sequence of internal connected sums of the $H_{g_i}$  together to form a manifold, say $H$,  that represents the given partition.  Using partitions into at most pairs and singletons, it follows from the $n=2$ case that rescaled tensor products of the form $$h^{-\lfloor \frac 1 2  \sum_i \varepsilon(\bb_i)  \rfloor }
\bb_1\otimes \cdots \otimes \bb_n \ \ \ \ \ \ (\bb_i\in \BB_i)$$
lie in $\BS(\sqcup_i \Si_i)$. Thus we only have to show that other partitions do not give rise to elements not in the $\BO$-span of these ones.
By induction on $n$, it is enough to consider the case where $H$ is connected. 
For this, we apply the same strategy as used in the proof for $n=2$.

Let $G_i$ be lollipop trees in $H_{g_i}$. We identify the  boundary connected sum of all the $H_{g_i}$ with $H_{g}$ where $g= \sum g_i$. The boundary of $H_{g}$ is $\#_i \Si_i$, the connected sum of all the  $\Si_i.$  Our $H$ is the interior connected sum $\# H_g$ with boundary $\sqcup_i \Si_i$; it is obtained from $H_g$ by adding $n-1$ two-handles  along the $n-1$ connecting circles in $\#_i \Si_i$ . 
 \begin{figure}[h]
\includegraphics[width=1.5in]{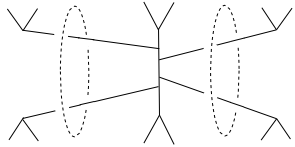}
\caption{Grafting  three  lollipop trees. If $n>3$, one generalizes this  pattern in the obvious way.} \label{3graft}
\end{figure}

We graft the $G_i$ together as indicated in Figure \ref{3graft} to form a lollipop tree $G$ in $H_{g}$. As in the case $n=2$, 
if $\bb$ is a basis element of $\BS(\#_i \Si_i)$ given by a small coloring of $G$, we denote by $\hat \bb$ the induced element of $\BS(\sqcup_i \Si_i)$ {\em via} the inclusion $H_g \subset H$. Our task is to compute the space spanned by $v$-graphs in $H$, but since they all can be isotoped into $H_g$, we only have to compute the span of  the elements 
$\hat \bb$ as  $\bb$ ranges over  the basis associated to the small colorings of $G$.

Again as in the case $n=2$, $\hat \bb$ will be zero unless $n-1$ specified pairs of edges have the same color. If these are the same, then we may define a related small coloring on each $G_i$. Let $\bb_i$ denote the
basis vector for $\BS(\Si_i)$ associated to this coloring.   As above,
  if   $\hat \bb$ does not  represent zero, 
 $\hat \bb $ is, up to units of $\BO$,  given by 
$ h^E \bb_1\otimes \ldots \otimes\bb_n$, where

 \begin{align*}
&E=   (n-1)(d-1)  -  \lfloor \frac 1 2 ( A(\bb)  -e(\bb)) \rfloor + \sum_i \lfloor 
\frac 1 2 ( A(\bb_i)  -e(\bb_i) ) \rfloor \\  
&  \ge   (n-1)(d-1)   -   \frac 1 2 ( A(\bb)  -e(\bb))  -   \frac n 2   +
             \sum_i  \frac 1 2 ( A(\bb_i)  -e(\bb_i) )   \\   
        &  = (n-1)(d-1)  -  \frac n 2 + \frac 1 2  ( e(\bb) -  \sum_i    e(\bb_i) ) 
         \text{,  as $A(\bb)=\sum_i  A(\bb_i)$} \\
        & \ge    (n-1)(d-1) -    \frac n 2  - \frac {n (d-1)} 2  
        \text{, as $e(\bb) \ge 0$, and   $e(\bb_i) \le  d-1
$} \\
       & =  \frac {(n-2)(d-2) } {2} -1.
       \end{align*} 
If $p >5$, then $d>2,$  and $E$ is non-negative (since $n\geq 3$ here).  Thus $v$-graphs in $H$ do not give rise to any new elements of  
$\BS(\sqcup_i \Si_i)$. 

We are reduced to the case that $p=5$
and $d=2$,
 when $E$ perhaps can be negative but is never less than $-1$. One has $E=-1$ exactly when
all $\geq$ signs in the computation above are equalities. This requires in particular that $A(\bb_i) -e(\bb_i)$ must be odd for all $i$, and that $e(\bb_i) =d-1$ for all $i$. In other words, we must have $\varepsilon(\bb_i)=1$ for all $i$.

{}From 
this we conclude that  for $p=5$ the connect-summing of $n\geq 3$ handlebodies together may imply the divisibility of an element of the form $\bb_1\otimes \cdots \otimes \bb_n$  by $h$ 
but in every case, one can also see at least this  divisibility arising from a partitioning of these handlebodies into at most pairs and singletons. This completes the proof.
 \end{proof}

\begin{rem}{\em It might have occurred to the reader that perhaps we should have defined $\BS$ of a disconnected surface by a tensor product formula, rather than as in Definition \ref{defS}.  This would have definite drawbacks when we consider how cobordisms act in the $\BS$-theory. Here is an example.

Suppose we let $\Si= \partial H_2$ with a point colored $p-3$ and let $C$ be the cobordism from $\Si \# \Si$ to $\Si \sqcup \Si$ constructed by adding a 2-handle as in Example \ref{ex}. The TQFT associates to $C$ a linear map $Z_p(C)$ from $V_p(\Si \# \Si)$ to  $V_p(\Si \sqcup \Si ).$  If we had defined $\BS(\Si \sqcup \Si )$ to be $\BS(\Si) \otimes \BS(\Si)$, 
then  $Z_p(C)$ would not map  $\BS(\Si \# \Si )$ into  $\BS(\Si \sqcup \Si ).$ 

On the other hand, with Definition \ref{defS} of $\BS$, for any cobordism $C$ from, say $\Si'$ to $\Si''$ such that $\beta_0(C,\Si'')=0,$ the induced map $Z_p(C)$ will  map 
 $\BS(\Si ')$ into  $\BS(\Si'').$ Such cobordisms are called {\em targeted} in \cite{G}.
More generally, we have that if $C$ is any cobordism $C$ from $\Si'$ to $\Si''$, then $Z_p(C)$  maps $\BS(\Si ')$ into  $\BD^{-\beta_0(C,\Si'')} \BS(\Si'').$

We remark that our theory is not half-projective with parameter $\BD$ in the sense of Kerler \cite{K}. Nor does it seem possible to make it so by rescaling.} \end{rem}

 \section{Extra 
structure, mapping
class group representations, and $\BSplus(\Si)$}\label{mcg.sec}
 
 The TQFTs we have been studying require that surfaces and $3$-manifolds be equipped with some extra structure in order to avoid having a framing anomaly. 
 Until now this extra structure did not play any essential role in our arguments. Thus a detailed discussion was not required.
 However some of the results in the following sections are improved by using the refined theory $V^+_p$ and its integral version $\BSplus$ which are defined and discussed in \cite{G,GMW}. For this reason, we now discuss the extra structures which we employ.
For simplicity, we will focus on the case where $\Si$ is connected in this and the next section (but we allow $\Si$ to have colored points).  
   
As in \cite{G,QG}, we 
follow the method originally developed by  Walker \cite{W} and Turaev \cite{T}. This approach equips $3$-manifolds with integer weights, and surfaces with Lagrangian subspaces of their first 
real 
homology. The cobordisms from $\Si$ to $\Si$ which are 
mapping cylinders 
form a central extension, denoted 
$\widetilde \Gamma(\Si)$,
 of the mapping class group. Forgetting the 
weight on these cobordisms 
defines a projection onto  the 
ordinary mapping class group $\Gamma(\Si).$ The kernel of this homomorphism is the group of integers $\BZ$. 

The extended mapping class group $\widetilde \Gamma(\Si)$ acts on $\BS(\Si)$ preserving the form $(\ ,\ )_{\Sigma}.$ Thus $\widetilde \Gamma(\Si)$ acts on $\BSs(\Si)$ as well. We note that elements of $\widetilde \Gamma(\Si)$ which map to the identity of $\Gamma(\Si)$ act on $\BS(\Si)$ by multiplication by some power of $\zeta_p$, 
 if $p \equiv -1 \pmod{4} $, and by some power of $\zeta_{4p}$, if  $p \equiv 1 \pmod{4}$.

By a standard handlebody, we mean one equipped with the weight $0$. The boundary of a standard handlebody $H$ will be equipped with the Lagrangian given by
the kernel of the inclusion to the homology of $H$.

In \cite[ \S 7]{G}, the notion of an even cobordism is defined. The even 
cobordisms from $\Si$ to $\Si$ which are 
mapping cylinders form an index two subgroup of  $\widetilde \Gamma(\Si)$ denoted  $\widetilde \Gamma(\Si)^+$ which still surjects onto $\Gamma(\Si)$.

In the remainder of this section, we assume $p \equiv 1 \pmod{4}$. Even cobordisms were used in \cite{G} to reduce coefficients from $\BO=\BZ[\zeta_{4p}]$ to $\BZ[\zeta_{p}]$ which we will denote by $\BO^+$. One obtains free $\BO^+$-modules $\BSplus(\Si)$ so that $$\BSplus(\Si)\otimes_{\BO^+}\BO=\BS(\Si)~.$$ The module $\BSplus(\Si)$ is again generated by $v$-graphs in a standard handlebody $H$, but coefficients are now required to lie in $\BO^+$. (Here we are using the fact that  a standard handlebody $H$ when viewed as a morphism from the empty set to $\Si$ is even.) Thus
our $\bb(a,b,c)$ also define elements of $\BSplus(\Si)$. Moreover the set of $\bb(a,b,c)$ associated to small colorings forms a basis for  $\BSplus(\Si)$ by the same proof. 

There is a sesquilinear form 
 \[ (\ ,\ )^+_{\Sigma}: \BSplus(\Sigma)
\times
\BSplus(\Sigma) \rightarrow \mathcal{O}^+\]
obtained by multiplying the form $(\ ,\ )_{\Sigma}$
by $i^{\delta (\Si)},$ where ${\delta (\Si)}$ is zero or one 
depending on whether the genus of $\Si$ is even, or odd. The reason why this form takes values in $\BO^+$ is explained in \cite[Remark 9.6]{GMW}. Note that $(\ ,\ )^+_{\Sigma}$ is hermitian or skew-hermitian, depending on the parity of the genus.
 
 We define $\BSps(\Si)$ in the same way as $\BSs(\Si)$ but using the form $(\ ,\ )^+_{\Sigma}$. Then the $\bb^\sharp(a,b,c)$ form a basis for $\BSps(\Si).$

The even extended mapping class group  $\widetilde \Gamma(\Si)^+$ acts on $\BSplus(\Si)$ preserving the form $(\ ,\ )^+_{\Sigma},$ and therefore also acts on $\BSps(\Si)$. Elements of $\widetilde \Gamma(\Si)^+$ which map to the identity of $\Gamma(\Si)$ act by multiplication by some power of $\zeta_p.$

  \section{Associated finite torsion modules}\label{torsion.sec}

The extended mapping class group $\widetilde \Gamma(\Si)$ acts on $\BS(\Si)/h^N \BS(\Si)$ ($N\geq 1$) and also on $\BSs(\Si)/\BS(\Si).$ These are finitely generated torsion $\BO$-modules. The hermitian form ${(\ , \ )}_\Si$ on $\BS(\Si)$ induces in the obvious way an $\BO/h^N \BO$-valued form on $\BS(\Si)/h^N \BS(\Si)$ and an 
$\BO(\frac 1 h) /\BO$-valued form on $\BSs(\Si)/\BS(\Si)$, and $\widetilde \Gamma(\Si)$  acts preserving these forms. The structure of these modules and forms follows easily from our bases. Let us describe them in some interesting cases. If $p \equiv 1 \pmod{4}$, we will mainly look at the modules coming from the refined theory $\BSplus$.

Let $G$ be a lollipop tree and consider the basis vectors $\bb(a,b,c)$ 
for $\BS(\Si)$, where $(a,b,c)$ runs through the small colorings of $G$.
 Clearly, $\BS(\Si)/h^N \BS(\Si)$ is a free 
$\BO/h^N \BO$-module with the induced basis.

\begin{prop} One has an orthogonal decomposition $$\BS(\Si) = \bigoplus_{a,c}
 \BS(\Si)_{G;a,c}$$ with respect to the form ${(\ , \ )}_\Si$, where $\BS(\Si)_{G;a,c}$ is the $\BO$-span  of the basis vectors $\bb(a,b,c)$ where $b$ varies so that $(a,b,c)$ is a small coloring of $G$.
\end{prop}
\begin{proof} This follows from the fact that the small graph basis vectors $\bg(a,b,c)$ are an orthogonal basis of $V_p(\Si)$ \cite[Theorem 4.11]{BHMV2}.
\end{proof}

 Clearly, this also induces an orthogonal decomposition of the induced form on $\BS(\Si)/h^N \BS(\Si)$.

Let us now look at $\BSs(\Si)/\BS(\Si)$. One has $\bb^\sharp(a,b,c)=h^{-n(a,c)}\bb(a,b,c)$, where $n(a,c)= e$, if $\sum a_i$ is even,  and $n(a,c)= e+1$, if $\sum a_i$ is odd (here, as always, $e$ is the trunk half-color of $(a,b,c)$). One has an orthogonal decomposition 
$$\BSs(\Si)/\BS(\Si)=
\bigoplus_{(a,c)}\BSs(\Si)_{G;a,c}/\BS(\Si)_{G;a,c}~,$$ and each summand is a free $\BO/h^{n(a,c)}\BO$-module.

 If $p \equiv 1 \pmod{4}$, similar statements hold for the $\BO^+$-modules 
$\BSplus(\Si)/h^N \BSplus(\Si)$ ($N\geq 1$) and  $\BSps(\Si)/\BSplus(\Si).$

We will say $(a,b,c)$ is an 
{\em odd (even)} small coloring if $\sum a_i$ is odd (even respectively).

If $p \equiv -1 \pmod{4},$ then $\BO=\BZ[\zeta_p]$ and  $\BO/h \BO = \BF_p$,  the finite field with $p$ elements.  If $p \equiv 1  \pmod{4}$,  then $\BO=\BZ[\zeta_{4p}]=\BZ[\zeta_{p},i]$ and $\BO/h \BO = \BF_p [i]$. 
In this case, we mainly consider $\BO^+/h \BO^+$ which is again equal to  $\BF_p$.

\begin{prop}\label{14.2} If $p \equiv -1 \pmod{4}$, the hermitian form ${(\ , \ )}_\Si$ induces a symmetric form on the $\BF_p$-vector space $\BS(\Si)/h\BS(\Si)$. If $p \equiv 1 \pmod{4}$, the $(-1)^g$-hermitian form $(\ , \ )^+_\Si$ induces a $(-1)^g$-symmetric form on the $\BF_p$-vector space $\BSplus(\Si)/h\BSplus(\Si)$. 
(Here $g$ is the genus of $\Si$.)
In both cases, one has an action of the (ordinary) mapping class group $\Gamma(\Si)$ preserving these forms.
\end{prop} 

\begin{proof} The symmetry properties of the induced forms follow from the fact that the conjugation on $\BZ[\zeta_p]$ induces  the trivial involution on $\BF_p$. The action of the (even) extended mapping class group descends to the ordinary mapping class group $\Gamma(\Si)$ because $\zeta_p$ acts as the identity on 
$\BF_p$.
\end{proof}

\begin{rem}\label{14.3}{\em If $\gamma \subset \Si$ is a separating simple closed curve, and all the colored points of $\Si$ lie on one side of $\Gamma$, then the associated Dehn twist $T_\gamma\in \Gamma(\Si)$ acts trivially  on the $\BF_p$-vector spaces above. This follows from Theorem~\ref{11.1}. Indeed,  $T_\gamma$ acts diagonally in the basis given by that theorem, and its eigenvalues are the twist
 coefficients $\mu_{2e}$ (see {\em e.g.} \cite[Remark 7.6(ii)]{BHMV2}) which are powers of $\zeta_p$, hence congruent to 
$1 \pmod{h}$. This argument is due to  Kerler.}\end{rem}

The forms given by Proposition~\ref{14.2} are singular in general. Of course 
the radical of such a form is preserved by the mapping class group. In this way we obtain, as well, a subrepresentation of characteristic $p$ given by the radical, and a quotient representation of characteristic $p.$ The quotient representation  preserves the induced  non-singular inner product space structure. We denote the radical of the form by $\rad_p(\Si).$

For the rest of this section, we consider only the important special case that  $\Si$ has no colored points. Then the radical is spanned by the images of the  $\bb(a,b,c)$ associated to odd small colorings. The quotient representation has dimension given by the number of even small colorings. 

Observe that $\BSs(\Si)/\BS(\Si)$  is also the free $\BO/h \BO$-module on the odd small colorings of $G.$

\begin{prop}\label{14.4} If $p \equiv -1 \pmod{4}$,  there is an isomorphism between the $\BF_p$-vector spaces
$\BSs(\Si)/\BS(\Si)$ and $\rad_p(\Si) $ intertwining $\Gamma(\Si)$ representations.  If $p \equiv 1 \pmod{4}$,  there is an isomorphism between 
$\BSps(\Si)/\BSplus (\Si)$ and $\rad_p(\Si) $ intertwining the $\Gamma(\Si)$ representations.
\end{prop}

\begin{proof} In the first case, consider the composition of  the map from $\BSs(\Si)$ to 
$\BS(\Si)$ given by multiplication by $h$ followed by the quotient map  to
$\BS(\Si)/h\BS(\Si).$ This is onto $\rad_p(\Si) $, and has  $\BS(\Si)$ as kernel. The second case is seen  similarly.
\end{proof}

Recall that the hermitian form ${(\ , \ )}_\Si$ on $\BS(\Si)$
induces a non-singular $\BO(\frac 1 h) /\BO$-valued form on $\BSs(\Si)/\BS(\Si)$. 
In our special situation, it suffices to multiply this form by $h$ to get an 
$\BO /h\BO$-valued form on $\BSs(\Si)/\BS(\Si)$. We denote this new form on $\BS(\Si)$ by $h.{(\ , \ )}_\Si$.

\begin{prop} If $p \equiv -1 \pmod{4}$, the form $h. {(\ , \ )}_\Si$ induces a skew-symmetric form on the $\BF_p$-vector space $\BSs(\Si)/\BS(\Si)$. If $p \equiv 1 \pmod{4}$, the form $h. (\ , \ )^+_\Si$ induces a $(-1)^{g+1}$-symmetric form on the $\BF_p$-vector space $\BSps(\Si)/\BSplus(\Si)$. In both cases, the forms are non-degenerate, and one has an action of the mapping class group $\Gamma(\Si)$ preserving these forms.
\end{prop}

\begin{proof} The proof is basically the same as for Proposition~\ref{14.2}. One just needs to observe that $ \bar h = 1- \zeta_p^{-1}= - \zeta_p^{-1} h $. Assume $p \equiv -1 \pmod{4}$. For $x,y\in \BSs(\Si)$ one has $$h.(y,x)_\Si= h. \overline{(x,y)_\Si}= - \zeta_p\, \overline{h.(x,y)_\Si}$$ as elements of $\BO$. As before, $\zeta_p$ acts trivially on $\BF_p$ and the induced conjugation is the identity. Only the minus sign remains. It is the reason why one gets a skew-symmetric form on $\BSs(\Si)/\BS(\Si)$ while the form in Proposition~\ref{14.2} was symmetric. The case $p \equiv 1 \pmod{4}$ is proved in the same way.
\end{proof}

Via the isomorphism of Proposition~\ref{14.4}, one gets in this way an inner product structure 
on $\rad_p(\Si) $ as well.

\section{An upper bound for the cut number of a $3$-manifold}\label{sec.cut}

The 
co-rank
of a group  is the maximal $k$ such that there is an epimorphism of that group onto the free group on $k$ letters \cite{S}.
If $M$ is an oriented connected closed 
$3$-manifold, define $c(M),$ the {\em cut number} of $M,$ to be the
maximal number of closed oriented surfaces that one can place in $M$ and still
have a connected complement.
Jaco \cite{Jaco} showed that the co-rank of $\pi_1(M)$ is the cut number 
of $M.$ It is easy to see that $c(M) \le  {\beta_1(M)}.$

In this section, we show that another upper bound on the cut number is given by quantum 
$SO(3)$-invariants.  Let $M$ be an oriented connected closed 
$3$-manifold and $L$ a banded colored graph in $M$. We use the normalization
$I_p(M,L) = \BD \, \langle (M,L)\rangle_p$ where $\langle \ \rangle_p$ is the invariant defined in \cite{BHMV2}, and 
$M$ is given the weight zero. This 
is
the same normalization as in \cite{MR}.
For example, $I_p(S^3)=1$ (since $\BD=\langle S^3\rangle_p^{-1}$). The invariant $I_p$ takes values in $\BO$ by \cite{Mu2,MR}.
We let $\bo_p(M,L)$  
denote
the highest power of $h$ which divides $I_p(M,L)$.

\begin{thm} \label{cut} Let $M$ be an oriented connected closed  
$3$-manifold,
$L$ a banded colored graph $L$ in $M$,  and 
$p= 2d+1 \ge 5$  
a prime.
 Then $$c (M) \le \frac {\bo_p(M, L)} {d-1}~.$$ 
\end{thm}

Varying $L$,  we obtain many upper bounds on the cut number of $M$. So far,  no example has been found where an interesting upper bound on $c (M)$ is obtained by choosing $L$ nonempty that cannot also be obtained by choosing $L$ empty.

Theorem \ref{cut} was conjectured  by the first author. He and  Kerler then obtained Theorem \ref{cut} for $p=5$  and  $L$ empty.  
 The proof we give of the more general result incorporates some ideas 
 that were
 used in 
 the
 proof of  
 this
 special case.

Cochran and Melvin showed \cite{CM} that 
$\beta_1(M)/3 \le  \bo_p(M)/(d-1)$.  
 Three recent papers
show that $\beta_1(M)/3$ can be larger than the cut number $c(M)$  \cite{H, LR,Si}. In this case the result of Cochran and Melvin gives a better lower bound on  $\bo_p(M)$  than Theorem \ref{cut}.  
On the other hand, when $\beta_1(M)/3 < c(M) $ (which may also happen), Theorem \ref{cut}  gives  a better lower bound on  $\bo_p(M)$.

In the proof, we will use the following easy proposition which follows from the functoriality properties of the TQFT. We will also use it in Section~\ref{sec.FKB}.

\begin{prop}\label{exp} Suppose $N \subset M$ and $N'$ is the exterior of $N$ in $M$. 
Further suppose $[N] \in  \BS(\partial N)$ can be written as a linear combination $\sum_i a_i [N_i]$.  We have that  
$I_p(M) =  \sum_i a_i  I_p(N_i \cup_{\partial N} N')$.
\end{prop}

\begin{proof}[Proof of Theorem \ref{cut}] 
Let $c$ be the cut number of $M.$
We can find $c$ disjoint  connected oriented surfaces in $M$ which  do not disconnect $M.$ As these surfaces do not disconnect, we may tube up parallel copies of these surfaces to find  a copy of the connected sum of these surfaces  disjoint from the original surfaces but homologous to the union of  these surfaces. If we add this new surface to the collection, we obtain  a collection  of $c+1$ connected oriented surfaces which disconnect but no sub-collection of $c$ of them disconnect. These $c+1$ surfaces disconnect $M$ into two components $Y$ and $Y'$ and each of the
$c+1$ surfaces is a boundary component of the closures of both $Y$ and $Y'$.

Let $\Si_i$ for $1 \le i \le c +1$ denote the individual surfaces.   We can find   $c$ disjoint arcs $ \alpha_i$  in $Y$ for  $1 \le i \le c$   which join  points $x_i \in \Si_i$ to points  $y_i \in \Si_{i+1}$.  Similarly we can find   $c$ disjoint arcs $ \alpha'_i$ in $Y'$ which join  
the points $x_i\in \Si_i$ to the points  $y_i \in \Si_{i+1}$.
A neighborhood of the union of the $\Si_i$ and the $\alpha_i$ and the $\alpha'_i$ is a cobordism $P$  from the connected sum of the $\Si_i$ to itself. We  isotope $L$ so that it is transverse to the $\Si_i$'s. We assign colored points to each $\Si_i$ according to this intersection.

We denote the connected sum of the $\Si_i$ by $\Si.$ On $\Si$ we have $c$ simple closed curves $\gamma_i$
such that if we perform surgery along them we recover the disjoint union $\sqcup_i \Si_i.$ The cobordism $P$ has a nice handle decomposition as follows:  $\Si \times I$ union $c$ 2-handles along $\{ \gamma_i \}_{ 1 \le i \le c}$   union $c$ 1-handles which ``reconnect''
the $\Si_i$.  The core of each 1-handle can be completed to a circle which meets exactly one core of a 2-handle transversely in a singe point.  We can identify $P$ with the result of framed surgery  
on
$\Si \times I$ along $\{ \gamma_i \times \frac 1 2 \}_{ 1 \le i \le c}$ with the framing along $\Si \times \frac 1 2$. Let $\gamma(\omega)$ denote $\{ \gamma_i \times \frac 1 2 \}_{ 1 \le i \le c}$ colored with the skein element $\omega$ appearing in the surgery axiom (S2) of \cite{BHMV2}. The surgery axiom says that $(\Si \times I, \gamma(\omega))$  induces the same endomorphism of $V_p(\Si)$ under the TQFT as does $P.$ 

Place $\Si \times I$ in $S^3$ so that its complement is the disjoint union of two handlebodies  $H$ and $H'$ and the $ \gamma_j$ bound disks $D_j$ in $H$. Express $[Y] \in \BS (\Si)$ in terms of the basis $\BB$ associated to a lollipop tree $G$ for $H$ with respect to $\ell(\Si)$. \footnote{Warning: we cannot assume that $G$ meets each $D_j$ in a single point, as there may be colored points on each $\Si_i$. But this causes no problem in the argument.} Similarly express $[Y'] \in \BS (-\Si)$ in terms of the basis $\BB'$ associated to a lollipop tree $G'$ for $H'$ with respect to $\ell(\Si)$.

Applying Proposition \ref{exp} twice,  once  to expand $Y$ and once to expand $Y'$, we have that $I_p(M,L)$ is a linear combination over $\BO$ of evaluations of skein classes  in $S^3$. Recall that every basis element $\bb=\bb(a,b,c)\in \BB$  is given by $ h^{- \lfloor\frac 1 2 ( A(\bb)-e(\bb))\rfloor} \bg(a,0,c)$ union some $v$-colored curves. Thus the term corresponding to $(\bb,\bb')\in \BB\times \BB'$ in the expansion of $I_p(M,L)$ is an integral multiple of   $ h^{- \lfloor\frac 1 2 ( A(\bb)-e(\bb) ) \rfloor - \lfloor\frac 1 2 (A(\bb')-e(\bb')) \rfloor }$ times the evaluation of the union of 

(i) a colored trivalent graph, namely the gluing of the graphs $\bg(a,0,c)$ in $H$ and $\bg(a',0,c')$ in $H'$ along their univalent vertices (which correspond to $\ell(\Si)$),

(ii) the $ \gamma_i$'s each colored by $\omega$, and 

(iii) some some extra 
$v$-colored
curves.

 We can do fusion along the strands which pass through the spanning disks $D_j$ for the $\omega$-colored $ \gamma_j$'s (none of the $v$-colored curves do), and discard the terms arising with a single non-zero color passing through these 
 disks, as for all these terms the color must be even.
 Then  each  $\omega$-colored $ \gamma_j$ spanning a disk may be replaced by 
an
extra scalar factor of $\BD$. (This is the same argument as in \cite[Proof of Lemma 4.1]{MR}.) Up to units this gives an extra factor of $h^{(d-1)c}.$ On the other hand, the Lollipop Lemma~\ref{ll}  
gives an extra factor of 
$h^{ \lceil\frac 1 2 (A(\bb) 
+ A(\bb')  )\rceil }$ which  compensates the above-mentioned negative power of $h$. Thus each term in the expansion of $I_p(M,L)$ is divisible by $h^{(d-1)c}.$ The  result follows. 
 \end{proof}
 
 \begin{ex}{\em Let $\Si$ be a closed surface of genus $2$, and let $T$ be the Dehn twist along an essential  separating curve 
$\gamma$.
Let $M_n$ be the mapping torus of $T^n$. Essential curves disjoint from $\gamma$ sweep out 
tori 
 in $M_n$. Picking one such curve on each side of $\gamma$, we see that  the cut number $c(M_n)\geq 2$. On the other hand, computing the trace of the map induced by $T^n$ on $V_p(\Si)$, we find $$I_p(M_n)=\BD \sum_{j=0}^{d-1} \zeta_p^{2nj(j+1)} (d-j)^2~.$$  If $n$ is not a multiple of $p$, this can be evaluated using Gauss sums. One finds that 
  $\bo_p(M_n)=2d-2$, and $c(M_n)\leq 2$ by Theorem~\ref{cut}. \footnote{If $n$ is a multiple of $p$, then $I_p(M_n)$ is $\BD$ times $\dim V_p(\Si)=d(d+1)(2d+1)/6$; this does not lead to a good upper bound for $c(M_n)$.}  Varying $p$, this shows $c(M_n)=2$ for every non-zero $n$. We remark that this can also be seen classically by computing the co-rank of $\pi_1(M_n)$.}\end{ex}

  \section{The Frohman Kania-Bartoszynska ideal invariant}
\label{sec.FKB}

\begin{de} [\cite{FK}] Given a connected $3$-manifold $N$ with boundary,
 let $\BJ_p(N)$  be the  ideal in $\BO$ generated by 
\[ 
\{ I_p(M)| \text{M is a closed connected oriented $3$-manifold containing $N$} \}. \]
In the case $p\equiv 1 \pmod{4},$ we also define
$\BJ_p^+(N)= \BJ_p(N)\cap \BO^+$. \end{de} 

\begin{rem}{ \em If $p\equiv 1 \pmod{4},$ the ideal $\BJ_p(N)$ is generated by scalars which are either in $\BO^+$ or in $i \BO^+$. (In fact, $I_p(M)$ lies in $i^{\beta_1(M)}\BO^+$ \cite[Remark 5.2]{MR}.) Thus $\BJ_p(N)$ is generated over $\BO$ by $\BJ_p^+(N)$. This is why we prefer to use $\BJ_p^+(N)$
if $p\equiv 1 \pmod{4}.$ 
}\end{rem}

This ideal is interesting because of the following immediate proposition.

\begin{prop} [\cite{FK}] \label{FK} If $N_1$ embeds in $N_2$,  then $\BJ_p(N_2) \subset \BJ_p(N_1),$
and $\BJ_p^+(N_2) \subset \BJ^+_p(N_1),$ if $p\equiv 1 \pmod{4}$. \end{prop}

\begin{rem}{ \em Frohman and Kania-Bartoszynska actually made this definition using  the $SU(2)$ theory in place of the $SO(3)$ theory used here.} \end{rem}

The ideal $\BJ_p(N)$ is hard to compute from its definition, because it involves the quantum invariants of infinitely many manifolds. 
Frohman and Kania-Bartoszynska were able to show that a related ideal (associated to the Turaev-Viro invariant at the third root of unity) is non-trivial 
\footnote{ We take the trivial ideal of a ring to be the ring itself.}
for the union of two solid tori glued together by identifying neighborhoods of $(2,1)$ curves on their boundary. 
But it seems that $\BJ_p(N)$ has never 
been computed exactly, except
 for manifolds $N$ with boundary a $2$-sphere.\footnote{However, sometimes it is possible to compute that the 
ideal is trivial or a power of $(h)$, 
without using Theorem~\ref{calc},
by making use of known theorems about the quantum invariants of  closed
3-manifolds. See \cite{G2} for examples of this.}

We can give a finite set of generators for $\BJ_p(N)$ using our bases.

\begin{thm}\label{calc} Let $N$ be an oriented connected $3$-manifold with boundary $\Si$, 
then $\BJ_p(N)$  is generated by the scalar products  $([N], \bb)_{\Si}$ as $\bb$ varies over a basis for $\BS(\Si).$ 
\end{thm}

This follows immediately from Proposition \ref{exp}. 
An example where $\BJ_p(N)$ is computed using Theorem~\ref{calc} will be given below.

\begin{rem}\label{algol} {\em In practice, if $N$ has connected boundary, we may present it as  surgery on a link in a handlebody $H$ standardly embedded in $S^3$.  Then we can just as well replace the form $([N], \bb)_{\Si}$ with the generalized Hopf pairing  $(([N], \bb'))_{H,H'}$ (see Section \ref{pB}) where  $\bb'$ are  elements of the graph-like basis for $\BS(-\Si)$ associated to  a lollipop tree in  the complementary handlebody $H'$. Thus $\BJ_p(N)$ can be computed from the evaluation of a finite number of explicitly given skein elements (consisting of $v$-graphs together with some $\omega$-colored curves) in $S^3$.}\end{rem}   

When one computes examples of  
$\BJ_p(N)$, 
one frequently gets the trivial ideal or some power of $(h)$. Here is an example where it is neither.
 \begin{figure}[h]
\includegraphics[width=1in]{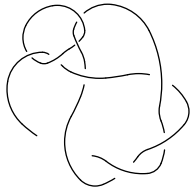}
\caption{L9a12 in 
Thistlethwaite's
list of prime links \cite{B} 
} \label{9a12} 
\end{figure}

 Let $L$ be the link  in  Figure \ref{9a12}. Let $K$ be the knotted 
component 
($K$ is a $(5,2)$ torus knot), 
and $J$ the unknotted component. Let $N(k)$ be the result of surgery to the exterior of $J$ along $K$ with framing $k$.
If $k$ is odd,  $N(k)$ is a homology circle.

 \begin{prop}\label{16.7} One has 
\begin{equation*}
\BJ^+_5(N(k))=\left\{
\begin{array}{cl}
(1 +2 \zeta_5^3)&\text{if $k \equiv 0 \pmod{5}~,$} \\
(1)&\text{otherwise}\\
\end{array}\right.
\end{equation*}
\end{prop}

\begin{proof}  
Proceeding as in  Remark \ref{algol}, the computation of the ideal
 $\BJ^+_5(N(k))$ in $\BZ[ \zeta_5]$ reduces to the evaluation of two skein elements in $S^3$
 (since $\Bs_5$ of a torus has rank two). We used data in Bar-Natan's Knot Atlas 
and his Mathematica package KnotTheory \cite{B} to help calculate this ideal.
\end{proof}

The ideal $(1 +2 \zeta_5^3)$ is 
a
non-trivial ideal as the norm of $1 +2 \zeta_5^3$ is $11$. Thus one has the following immediate

\begin{cor}\label{16.8} The manifold  $N(5n)$ does not embed into any closed $3$-manifold $M$  whose $I_5(M)$ is not divisible by 
 $1 +2 \zeta_5^3.$
 In particular $N(5n)$ does not embed into the 3-sphere 
(since $I_p(S^3)=1)$.
\end{cor}

Further such examples  are explored in \cite{G2}. Some more general results about the ideals $\BJ_p(N)$ are given there as well.

\end{document}